\documentclass[10pt]{amsart}
\usepackage{amsmath,amssymb}

\usepackage{latexsym}
\usepackage[dvips]{graphics}
\newtheorem{theorem}{Theorem}[section]
\newtheorem{proposition}{Proposition}[section]
\newtheorem{lemma}{Lemma}[section]
\newtheorem{corollary}[theorem]{Corollary}
\newtheorem{example}[theorem]{Example}

\theoremstyle{remark}
\newtheorem{remark}{Remark}[section]
\newtheorem{definition}{Definition}[section]

\def\p1i{\pi_1^{\infty}}

\let\ =\quad

\def \Z{\mathbb Z}
\def \R{\mathbb R}

\title{On the wgsc and qsf tameness conditions for finitely 
presented groups}

\author{Louis Funar}
\address{Institut Fourier BP 74, UMR 5582,  Universit\'e Grenoble I, 38402
Saint-Martin-d'H\`eres Cedex, France }
\email{funar@fourier.ujf-grenoble.fr}
\author{Daniele Ettore Otera}
\address{Universit\'a di Palermo, Dipartimento di Matematica e Applicazioni,
via Archirafi 34, 90123 Palermo}
\email{oterad@math.unipa.it}
\begin{document}


\begin{abstract}

A finitely presented group is weakly geometrically simply connected (wgsc) 
if it is the fundamental group  of some compact polyhedron  
whose universal covering is wgsc i.e. it  has 
an exhaustion by compact  connected and simply 
connected sub-polyhedra. 
We show that this condition is almost-equivalent to Brick's qsf property, which 
amounts to finding an exhaustion approximable by finite simply connected 
complexes, and also to the tame combability introduced and studied by 
Mihalik and Tschantz.  
We further observe that a number of standard constructions 
in group theory yield qsf groups and analyze specific 
examples. We show that requiring the exhaustion    
be made of metric balls in some Cayley complex 
is a strong constraint, not verified by general qsf groups.  
In the second part of this paper we give sufficient conditions  
under which  groups which are extensions 
of finitely presented groups by finitely generated 
(but infinitely presented) groups are qsf.  
We prove, in particular, that the finitely presented 
HNN extension of the Grigorchuk group is qsf.  

\vspace{0.1cm} \noindent {\bf Keywords:} Weak geometric simple
connectivity,  quasi-simple filtration,  tame combable, Grigorchuk group, HNN extension.

\vspace{0.1cm} \noindent {\bf MSC Subject:} 20 F 32, 57 M 50.
\end{abstract}

\maketitle

\section{Introduction}
\noindent  
Casson and Poenaru (\cite{Po2}, \cite{GS}) studied 
geometric conditions on the Cayley graph of a  finitely presented group 
implying that the universal covering of a compact 3-manifold 
with given fundamental group is $\R^3$.  The proof involves approximating 
the universal covering by compact, simply-connected three-manifolds.
This condition was then adapted for arbitrary spaces and finitely 
presented groups by S.Brick in \cite{BM1} (see also \cite{St2}), 
under the name {\em quasi-simply filtered} (abbreviated {\em qsf} below).

\vspace{0.2cm}
\noindent   
We consider  here a related and apparently
stronger notion, called {\em weak geometric simple connectivity} 
(abbreviated  {\em wgsc}), which came out from the 
study of  the geometric simple connectivity of open manifolds in \cite{FG}.
Specifically, a polyhedron is wgsc if it admits an exhaustion 
by compact connected and simply connected polyhedra.    
The interest of such a strengthening is that 
it is easier to prove that  specific high 
dimensional polyhedra are not wgsc rather than not qsf. 
In fact,  a major difficulty encountered 
when searching for examples of manifolds which are 
not wgsc  is that one has to show that no exhaustion 
has the required properties, while, in general, non-compact manifolds 
are precisely described by means of one specific exhaustion. 
Thus, one needs a method to decide 
whether a space is  not wgsc out of a given (arbitrary) exhaustion. 
We are not aware about such methods in the qsf setting.   
However, the criterion given in \cite{FG} permits  
to answer this question for the wgsc condition, 
at least for non-compact manifolds 
of high dimensions.

\vspace{0.2cm}
\noindent 
A central issue in geometric group theory is to study classes of groups 
with various properties of topological nature. 
The topological properties in question 
are imported from the realm of  infinite complexes by means of the following 
recipe, which was first used on a large scale by Gromov. 
Say that a finitely presented 
(in general infinite) group has a certain property 
if the universal covering of some finite complex with this 
fundamental group has the required property. In this setting 
we can speak about the qsf (or wgsc) of finitely presented groups. 
In this respect we have three levels of equivalence relations 
among topological properties. First, the usual one concerning 
(more or less) arbitrary CW complexes. Second, the {\em almost-equivalence} 
which concerns only universal coverings of finite complexes, i.e. 
finitely presented groups.  At last we have the quasi-isometry 
equivalence relation for finitely presented groups. 
In this paper we will mostly consider the almost-equivalence of various 
tameness properties, which will also permit us to draw conclusions 
about their quasi-isometry invariance. 
Notice however that qsf and wgsc have different flavors. 
If one universal covering of a finite complex with given fundamental 
group is qsf  then all such universal coverings are qsf and, 
in particular, this holds for every Cayley complex. 
Thus the qsf property is independent on the presentation used in 
the construction of the Cayley complex. 
This is not anymore true for the wgsc property. There are  
examples of presentations of a wgsc  group which lead to non 
wgsc Cayley complexes. However we
will see that these two properties define the same class of groups in the
sense that a group is qsf if and only if it is wgsc. The qsf is then a
group property which is presentation independent and almost-equivalent to
the wgsc.

\vspace{0.2cm}
\noindent 
The wgsc property should be compared to a tameness condition  
which is central in non-compact manifold theory, 
namely the {\em simple connectivity at infinity}. Roughly speaking 
the simple connectivity at infinity expresses the fact that loops which are 
far away  should bound disks which are far away. This topological 
property have been used for characterizing Euclidean spaces as being the 
contractible manifolds that are simply connected at infinity 
by Siebenmann, Stallings and Freedman. 
Moreover, the simple connectivity at infinity is much 
stronger than the wgsc in dimensions at least 4, and in particular 
for finitely presented groups.  In fact, 
M.Davis  (\cite{Da1}) constructed  examples of  
aspherical manifolds  whose  universal coverings are different 
from $\R^n$ (for $n\geq 4$). Further one understood 
that these examples are quite common (see \cite{DM}). 
The groups in these examples are finitely generated Coxeter groups, which act 
properly co-compactly on some CAT(0) complexes and thus they are wgsc.

\vspace{0.2cm}
\noindent 
In order to give an unified proof that many classes of groups are qsf 
Mihalik and Tschantz (\cite{MT2}) introduced the related notion 
of {\em tame 1-combings} for groups.  An usual combing for a 2-complex 
is the choice of paths in the 1-skeleton joining a base-point vertex to 
every other vertex.  Groups whose Cayley graphs admit nice (e.g. bounded) 
combings have good  algorithmic properties, like automatic groups and 
hyperbolic groups and were the subject of extensive study in the last 
twenty years. Further a 1-combing corresponds to one dimension higher, 
namely, to a system of paths  joining a base-point vertex to every point 
of the 1-skeleton. We refer to the next section for the precise definition   
(see 2.4) of the enhanced notion of tame 1-combing 
of 2-complexes (and groups). One of the main results of \cite{MT2} 
is that tame 1-combable groups (and in particular 
asynchronously automatic groups and semi-hyperbolic groups) are actually qsf.

\vspace{0.2cm}
\noindent 
Our aim is to pursue further the study of the qsf 
condition  for groups. 
The first part of this paper is devoted  
to finding  characterizations of the qsf 
by means of methods from high dimensional 
manifold theory. 
Our first result is the following.  

\begin{theorem}\label{main1}
The wgsc, qsf and tame 1-combability conditions are almost-equivalent 
topological properties of finitely presented groups. 
\end{theorem}

\vspace{0.1cm}
\noindent 
In particular, using  the results from \cite{Brick}, we obtain that: 

\begin{corollary}\label{qi} 
A group quasi-isometric to a qsf finitely presented group is qsf. 
\end{corollary}
  
\vspace{0.1cm}
\noindent 
In other words, the qsf property of groups is 
{\em geometric}. We apply these results  to analyze 
several interesting classes groups and derive additional examples  
of qsf groups. 

\vspace{0.2cm}
\noindent 
A natural question is whether there is some natural simply connected   
exhaustion  for a wgsc group. A possible candidate is to consider 
the word metric on the Cayley complex associated to a group 
presentation and the associated exhaustion by 
metric balls. We will show in section 4 that:   

\begin{theorem}\label{main2}
Finitely presented groups admitting a Cayley complex whose metric balls 
have fundamental groups generated by loops of uniformly bounded length 
have linear connectivity radius and solvable word problem. 
\end{theorem}

\vspace{0.cm}
\noindent  
In particular such groups are strongly constrained and there are  
examples of wgsc groups not satisfying these conditions.  
Therefore the simply connected exhaustions of wgsc Cayley complexes 
are far from being the ones  by metric balls. 
As application we will give a simple proof for the fact that  
finitely presented groups admitting complete geodesic rewriting systems 
are qsf.

\vspace{0.2cm}
\noindent 
The start point of section 5 
is the result of Brick and Mihalik from \cite{BM2} 
which states that extensions of  infinite finitely presented groups by  
finitely presented groups are qsf. This is the group theoretical analog 
of the fact that products of contractible manifolds are homeomorphic to the 
Euclidean space. Using the same methods we can prove the following:

\begin{theorem}\label{main4}
An ascending HNN extension of a finitely presented group 
is qsf. 
\end{theorem}

\vspace{0.2cm}
\noindent 
We investigate further extensions 
of infinite finitely presented groups by suitable  
infinitely presented groups, for instance torsion groups.  
The second main result of section 5 is Theorem \ref{ext} which 
gives  sufficient (too technical to state here) 
conditions for  such extensions to be qsf.  

\vspace{0.2cm}
\noindent 
Then we consider in detail the case of the 
Grigorchuk group and of its finitely 
presented  HNN extension constructed in 
(\cite{Gri0,Gri}, see also \cite{DLH}). 
The main result of the second part (see section 6) is that: 

\begin{theorem}\label{main3}
The finitely presented HNN extension of the Grigorchuk group 
is qsf. 
\end{theorem}
 
\vspace{0.1cm}
\noindent
These methods could be used 
in slightly more general situations in order to cover large classes 
of finitely presented extensions of branch groups having endomorphic 
presentations, as defined by Bartholdi in \cite{Bar}. However, the present 
approach does not permit to prove the qsf of all such extensions, without 
an additional condition.

\vspace{0.2cm}
\noindent 
At this point we wish to emphasize 
the difference between the geometric invariants of discrete groups 
and those of topological nature.  
Geometric invariants are sensitive to cut and paste 
operations and thus algebraic constructions  can provide 
a large variety of examples. For instance the set 
of exponents of polynomial isoperimetric Dehn functions 
of finitely presented groups is a dense subset of $[2,\infty)$. 
These correspond to distinct quasi-isometry classes 
of groups.   
On the other side,  topological properties are quite stable and thus       
can be satisfied by very large classes of groups.  Two typical cases 
are the semi-stability at infinity (see e.g. \cite{Mi3}) and 
the property that $H^2(G,\Z G)$ is free abelian. 
It is still unknown whether all finitely presented groups
satisfy either one of these two properties.

\vspace{0.2cm}
\noindent 
In the same spirit there are still no known examples of 
finitely presented groups which are not qsf (see \cite{St2}). 
Notice that fundamental groups of compact 3-manifolds  
are qsf, but the only proof of that is as a consequence 
of the Thurston geometrization conjecture (settled by Perelman).
If non qsf groups do exist they would lay  at the opposite 
extreme to hyperbolic and non-positively 
curved groups and thus  they should be highly non generic.  
A related question is whether   
the fundamental group of a closed aspherical manifold  could 
act properly (not necessarily co-compactly) 
on some non-wgsc contractible manifold, like those described in \cite{FG}.

\vspace{0.2cm}
\noindent 
{\bf Acknowledgements.} The authors are thankful to 
Tullio Ceccherini-Silberstein, Rostislav Grigorchuk, 
Francisco Lasheras, Pierre Pansu, Valentin Poenaru, 
Dusan Repovs, Tim Riley, Mark Sapir and the referee  
for useful comments. We are greatly indebted to Mike Mihalik for his 
careful reading of previous versions of the manuscript 
leading to numerous suggestions and corrections. 
L.F. was partially supported by 
Proteus program (2005-2006), no 08677YJ, the 
ANR Repsurf: ANR-06-BLAN-311 and D.O. by
GNSAGA of INDAM, and both authors by the Projet
``Internazionalizzazione''
\emph{Propriet\`a asintotiche di variet\`a e di gruppi
discreti} of Universit\`a Degli Studi di Palermo and MIUR of Italy.

\section{Preliminaries on tameness conditions for groups}

\subsection{The wgsc.}
\noindent The following definition due to  C.T.C.Wall came out
from the work of S.Smale on the Poincar\'e Conjecture and, more
recently, in the work of V.Poenaru (\cite{Po2}). Moreover, it has
been revealed as especially interesting in the non-compact
situation, in connection with uniformization problems (see
\cite{FG}).

\begin{definition}
A non-compact manifold, which might have nonempty boundary, is
{\em geometrically simply connected}  (abbreviated  {\em gsc})  if
it admits a proper handlebody decomposition without 1-handles, 
or equivalently, in which every 1-handle is in cancelling position with a 
2-handle. Alternatively, there 
exists a proper Morse function $f:W\to \R$, whose critical points 
are contained within  ${\rm int}(W)$ such that: 
\begin{enumerate}
\item  $f$ has no index one critical points; and  
\item the restriction $f|_{\partial W}:\partial W\to \R$ 
is still a proper Morse function without non-fake index one 
critical points. The non-fake critical points of  $f|_{\partial W}$ 
are those for which the gradient  vector field ${\rm grad}\, f$ points towards 
the interior of $W$, while the fake ones are those for which 
${\rm grad}\, f$ points outwards.  
\end{enumerate}   
\end{definition}

\noindent 
The gsc condition was shown to be a powerful tameness condition 
for open three--manifolds and four--manifolds in the 
series of papers by Poenaru starting with \cite{Po2}.

\begin{remark}
Handle decompositions are known to exist for all manifolds in the
topological, PL and smooth settings, except in the case of
non-smoothable topological 4-manifolds. Notice that open
4-manifolds are smoothable.
\end{remark}

\noindent Manifolds and handlebodies considered below are PL.

\vspace{0.2cm} \noindent  One has the following combinatorial
analog of the gsc for polyhedra:

\begin{definition}
A non-compact polyhedron $P$ is  {\em weakly geometrically simply
connected}  (abbreviated  {\em wgsc}) if $P=\cup_{j=1}^{\infty} K_j$,
where $K_1\subset K_2\subset \cdots\subset K_j\subset\cdots$ is an
exhaustion by compact connected sub-polyhedra with $\pi_1(K_j)=0$.
Alternatively, any compact sub-polyhedron is contained in a simply
connected sub-polyhedron.
\end{definition}

\vspace{0.2cm} \noindent Notice that a wgsc polyhedron is simply connected.  
The wgsc notion is the counterpart in the polyhedral category 
of the gsc of open manifolds and in general it is slightly weaker. 
The notion which seems to capture the full power of the gsc 
for non-compact manifolds (with boundary) is  the pl-gsc discussed 
in \cite{FLR}. 

\begin{remark}
Similar definitions can be given in the case of
topological (respectively smooth) manifolds where we require the
exhaustions to be by topological (respectively smooth)
sub-manifolds.
\end{remark}

\begin{remark}
For $n\neq 4$ an open $n$-manifold is wgsc if and only if it is gsc 
(see \cite{FG} for $n\geq 5$, and for $n=3$ 
it follows from the Poincar\'e conjecture). While in dimension 4 one
expects to find  open 4-manifolds which are wgsc but not gsc.  
\end{remark}

\begin{definition}
The finitely presented group $\Gamma$ is {\em wgsc} if {\em there exists
some} compact polyhedron $X$ with $\pi_1(X)=\Gamma$ such that
its universal covering $\widetilde X$  is wgsc.
\end{definition}

\begin{remark}
Working with simplicial complexes instead of polyhedra in the definitions
above, and thus not allowing subdivisions, 
yields an equivalent notion of wgsc for finitely presented groups.
\end{remark}

\begin{remark}
The fact that a group is not wgsc cannot be read from an arbitrary 
complex with the given fundamental group.
In fact, as F.Lasheras pointed out to us, 
for any  finitely presented group $\Gamma$ with an element of
infinite order, there exists a complex $X$ with $\pi_1 X=\Gamma$
whose  universal covering  is not wgsc. For example take $\Gamma=\Z$ 
and the complex $X$ being that associated to 
the presentation $\Z=\langle a,b| baba^{-1}b^{-1}\rangle$. Then the 
universal covering  $\widetilde{X}$ is: 

\begin{center}
$\includegraphics{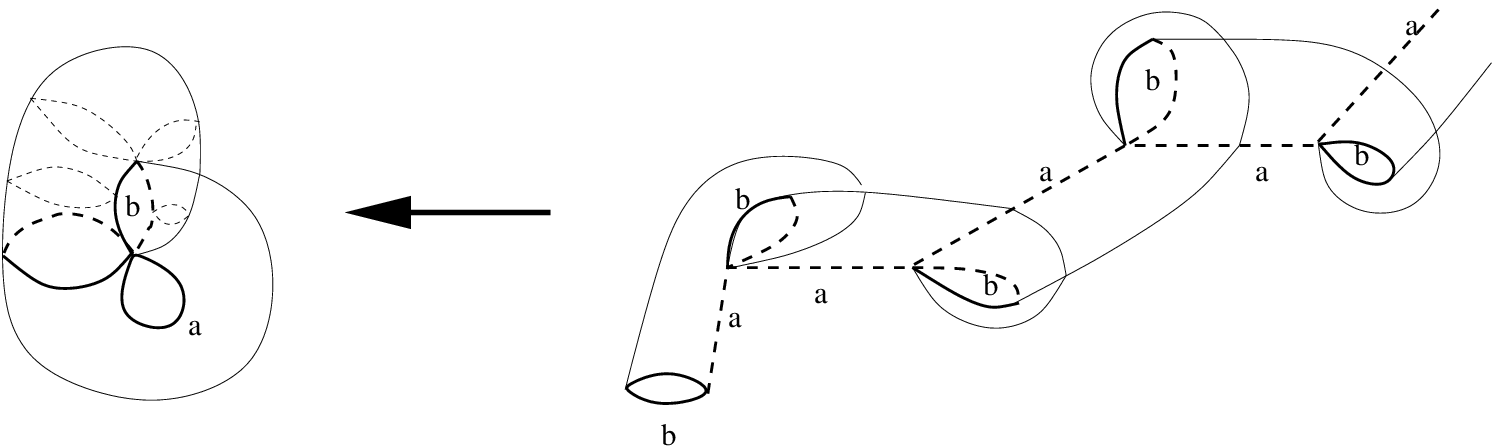}$
\end{center}
\noindent 
One sees that $\widetilde{X}$  is not wgsc because 
in the process of killing one loop $b$ one creates another one indefinitely.  

\vspace{0.2cm}
\noindent 
Further, if $\Gamma$ is a finitely presented group with an element $a$ 
of infinite order we add a new generator $b$  
and a new relation as before. The universal covering associated to this 
presentation is not wgsc, by the same arguments. 
\end{remark}

\begin{remark}
The wgsc property cannot be extended to arbitrary 
finitely generated groups, as stated, 
since any group admits a presentation with
infinitely many relations such that the associated 2-complex 
is wgsc. It suffices to add infinitely
many $2$-cells, along the boundaries of unions of 
$2$-cells, killing inductively the fundamental group of any compact subset.
\end{remark}

\begin{remark}
Recall that there exist uncountably many  open 
contractible manifolds which are not wgsc (\cite {FG}). In general, 
these manifolds are not covering spaces and we don't know whether 
one could find co-compact universal coverings among the non wgsc manifolds.
For instance, if a finitely presented torsion group exists then 
it is hard to believe that its Cayley complex is wgsc. 
Swenson has shown that every CAT(0) group has an element of 
infinite order (see \cite{Swe}).   
Notice that manifolds that are simply connected at infinity are 
automatically wgsc (\cite{PT}), but in general not conversely (see examples 
below). 
\end{remark}

\subsection{The qsf property after Brick and Mihalik}

\noindent The qsf property is a weaker version of the wgsc, which
has the advantage to be independent on the polyhedron we chose.
Specifically, Brick  (\cite{BM1}) defined it as follows:

\begin{definition}
The simply connected non-compact PL space $X$ is {\em qsf} if for
any compact sub-polyhedron $C\subset X$ there exists a simply connected compact
polyhedron $K$ and a PL  map $f: K\to X$ so that $C\subset f(K)$
and $f|_{f^{-1}(C)}:f^{-1}(C)\to C$ is a PL homeomorphism.
\end{definition}

\begin{definition}
The finitely presented group $\Gamma$ is {\em qsf} if there exists
a compact polyhedron $P$ of fundamental group $\Gamma$ so that its
universal covering $\widetilde{P}$ is qsf.
\end{definition}

\begin{remark}
 It is known (see \cite{BM1}) that the qsf  is a group property and
does not depend on the compact polyhedron $P$ we chose in the
definition above. In fact, if $Q$ is any compact polyhedron of
fundamental group $\Gamma$ (which is qsf) then $\widetilde{Q}$ is qsf.
\end{remark}

\begin{remark}
The qsf is very close to (and a consequence of)  
the following notion of Dehn exhaustibility (see \cite{Po2,FG}) 
which was mainly used in a manifold setting. 
The polyhedron $W$ is Dehn-exhaustible if for any compact 
$C\subset W$ there exists a simply connected compact polyhedron 
$K$ and an {\em immersion} $f:K\to W$  such that 
$C \subset f(K)$ and the  set of double points 
of $f$ is disjoint from $C$. It is known from \cite{Po2} that 
a Dehn exhaustible 3-manifold is wgsc.  
\end{remark}

\subsection{Small content and 1-tame groups}
Now we consider some other tameness conditions on non-compact
spaces, which are closely related to the wgsc. Moreover we will
show later that they induce equivalent notions for discrete
groups. In many cases it is easier to prove that a specific complex 
has one of these two properties instead that directly proving the qsf.
This will be the case in the second part of this paper for the 
Grigorchuk group and its extension.

\begin{definition}\label{small}
The simply connected non-compact polyhedron $X$ has {\em small
content} if for any compact $C\subset X$ there exist two compact
connected sub-polyhedra $C\subset D\subset E\subset X$, fulfilling the
following properties:
\begin{enumerate}
\item The map $\pi_1(D)\to \pi_1(E)$ induced by the inclusion, is
  zero.
\item If two points of $D$  are connected within $E-C$ then they are 
connected within $D-C$. 
\item Any loop in $E-C$ (based to a point in $D-C$)
is homotopic  rel the base point within $X-C$ to a loop which lies
entirely inside $D-C$. Alternatively,  let us denote by
$\iota_Y:\pi_1(Y-C)\to \pi_1(X-C)$ the morphism induced by
inclusion (for any compact $Y$ containing $D$), by fixing a base
point (which is considered to be in $D-C$). Then one requires that
$\iota_D(\pi_1(D-C))= \iota_E(\pi_1(E-C))$.
\end{enumerate}
The finitely presented group $\Gamma$ has {\em small content} if
there exists a compact polyhedron $P$ of fundamental group
$\Gamma$ so that its universal covering $\widetilde{P}$ has small
content.
\end{definition}

\begin{remark}
An obvious variation would be to ask that the homotopy above might
not keep fixed the base point. We don't know whether the new
definition is equivalent to the former one.
\end{remark}

\begin{definition}
The PL space $X$ is  {\em 1-tame} if any compact sub-polyhedron $C$
is contained in a compact sub-polyhedron $K\subset X$, so that any loop $\gamma$
in $K$ is (freely) homotopic within $K$ to a loop $\overline{\gamma}$ in $K-C$,
while $\overline{\gamma}$  is null-homotopic within $X-C$.

\vspace{0.1cm}
\noindent
The finitely presented group $\Gamma$ is {\em 1-tame} if
there exists a compact polyhedron $P$ of fundamental group
$\Gamma$ so that its universal covering $\widetilde{P}$ is 1-tame.
\end{definition}
\noindent Notice that one does not require that an arbitrary loop in $K-C$
be null-homotopic within $X-C$. This happens only after a suitable
homotopy which takes place in $K$.


\subsection{Tame combings and the Tucker property}
Group combings were essential ingredients in Thurston's attempt to 
abstract finiteness properties of fundamental groups of negatively curved  
manifolds which finally led to automatic groups.

\vspace{0.2cm}
\noindent 
Tame 1-combings of groups were considered by Mihalik and Tschantz in \cite{MT2} 
as higher dimensional analogs of usual combings, which are referred of 
as 0-combings. 

\begin{definition}
A $0$-{\em combing} of a 2-complex $X$ is a set of  
edge-paths $\sigma_p(t) $, $t\in[0,1]$, joining 
each vertex $p$ of $X$ to a 
base-point vertex $x_0$. 
This can be thought of as a homotopy 
$\sigma\colon X^0 \times [0,1] \rightarrow X^1$ for which 
$\sigma(x,1)=x$ for all $x \in X^0$, and 
$\sigma(X^0,0)=x_0$, where $X^j$ denotes the $j$-dimensional skeleton of 
$X$.  

\vspace{0.1cm}\noindent 
A $1$-{\em combing} of the 2-complex $X$ is a 
continuous family of paths $\sigma_p(t)$, $t\in[0,1]$,   
joining each point  $p$ of the 1-skeleton of $X$ to a base-point 
vertex $x_0$, whose restriction to vertices is a $0$-combing. 
This is a homotopy  
$\sigma\colon X^1 \times [0,1] \rightarrow X$ for 
which $\sigma(x,1)=x$ for all 
$x \in X^1$, $\sigma(X^1, 0)=x_0$, and 
$\sigma |_{X^0 \times [0,1]}$ is a $0$-combing. 
\end{definition}

\noindent Observe that although any connected complex is $0$-combable,  
a 2-complex is 1-combable if and only if it is simply connected. 

\vspace{0.2cm}
\noindent 
In order to find interesting consequences in geometric group theory 
one imposed the boundedness (or fellow traveler condition) on    
the $0$-combing, namely that combing paths of 
neighbor vertices be at uniformly bounded distance from each other. 

\vspace{0.2cm}
\noindent 
In the same spirit Mihalik and Tschantz replaced the boundedness by  
the following property of topological nature:

\begin{definition}
A $0$-combing is called {\em tame} if for every compact set 
$C \subseteq X$ there exists a compact set 
$K \subseteq X$ such that for each $x \in X^0$ 
the set $\sigma^{-1}(C) \cap (\{x\} \times [0,1])$ is contained in one path 
component of $\sigma^{-1}(K) \cap (\{x\} \times [0,1])$. 

\vspace{0.1cm}
\noindent A $1$-combing is {\em tame}  if its restriction to the  
set of vertices is a tame $0$-combing and for each compact 
$C\subset X$ there exists a larger compact $K\subset X$ such that 
for each edge $e$ of $X$, $\sigma^{-1}(C) \cap (e \times [0,1])$ 
is contained in one path component of $\sigma^{-1}(K) \cap (e \times [0,1])$.

\vspace{0.1cm}
\noindent A group is {\em tame 1-combable} if the universal cover of some 
(equivalently any, see \cite{MT2}) 
finite 2-complex with given fundamental group admits a tame 1-combing. 
\end{definition}

\vspace{0.1cm} 
\noindent Recall now the following tameness condition of topological 
spaces:

\begin{definition}
The non-compact PL space $X$ is {\em Tucker} if the fundamental group 
of each component of $X-K$ is 
finitely generated, for any finite sub-complex $K\subset X$. 
\end{definition}

\vspace{0.1cm}
\noindent 
This definition was motivated by Tucker's work \cite{Tu} on 3-manifolds. A 
non-compact manifold is a {\em missing boundary} manifold  if it is obtained 
from a compact manifold with boundary by removing a closed subset 
of its boundary.  We have the following characterization from \cite{Tu}:  
a $P^2$-irreducible connected 3-manifold is a  missing boundary 3-manifold
if and only if it is Tucker.

\noindent The main results of \cite{MT2} state that: 
\begin{proposition}[\cite{MT2}] \label{tamecomb}
A finitely presented group is tame 1-combable if and only if the 
universal covering of any (equivalently, some) finite complex 
with given fundamental group is Tucker. 
Moreover,  a tame 1-combable group is qsf. 
\end{proposition}

\noindent All known examples of qsf groups are actually tame 1-combable. 
We will show in the next section that the two notions are almost-equivalent. 

\vspace{0.1cm}
\noindent  Requiring a tame  $0$-combing is a very soft condition, since: 

\begin{proposition}
Any connected 2-complex $X$ has a tame $0$-combing. 
\end{proposition}
\begin{proof}
The key-point is that  any connected 2-complex $X$ is the ascending union 
of {\em connected} finite sub-complexes $X_n$, for instance metric balls.  

\vspace{0.1cm}
\noindent
A $0$-combing $\sigma_p$ is geodesic with respect to 
$(X_n)_n$ when it satisfies the following properties: 
\begin{enumerate}
\item if $p\in X_0$ then $\sigma_p$ has minimal length among 
the paths in $X_0$ joining $p$ to $x_0\in X_0$; 
\item for $p\in X_n-X_{n-1}$, with $n\geq 1$, there is some $q\in X_{n-1}$  
which realizes the distance in $X_n$ from $p$ to $X_{n-1}$.  Let $\eta_p\subset X_n$ 
be a minimal length curve joining $p$ to $q$. Then $\sigma_p$ is the concatenation  of 
$\eta_p$ and $\sigma_q$.  
\end{enumerate}

\vspace{0.1cm}
\noindent
If all $X_n$ are connected then there exist geodesic 
$0$-combings  which are defined inductively by means of the two conditions above. 
Let $\sigma$ be one of them. 
It suffices to verify the tameness of $\sigma$ for large enough 
finite sub-complexes  $C$ and thus to assume that $x_0\in C$. 
Set  $K$ be the smallest $X_n$ containing $C$. 

\vspace{0.1cm}
\noindent
We claim that the set $\{t\in [0,1]; \sigma_p(t)\in K\}$ is connected, 
 which settles the proposition.
This is clear when $p\in X_n$. If $p\in X_{n+1}-X_n$ then $\eta_p$ is contained in 
$X_{n+1}-X_n$ except for its endpoint $q\in X_n$. Otherwise, we would find a point 
in $X_n$ closest to $p$ than $q$ contradicting our choice for $q$.  
Further $\sigma_q\subset X_n$ is connected, hence  $\sigma_p\cap X_n$ is also connected. 
Using induction on $k$ one shows  in the same way that    
$\sigma_p\cap X_n$ is connected when $p\in X_{n+k}-X_n$. 
As $X=\cup_{k} X_{n+k}$ the claim follows.  
\end{proof}

\begin{remark}
One says that a 1-combing is weakly tame 1-combable if for each compact 
$C\subset X$ there exists a larger compact $D\subset X$ such that  for every edge $e$ the set 
$\{(p,t)\in e\times [0,1] ; \sigma_p(t)\in C\}$ is contained in 
one connected component of $\{(p,t)\in e\times [0,1]; \sigma_p(t)\in D\}$.  
Thus one drops from the definition of the tame 1-combing the requirement that 
the restriction to the vertices be a tame $0$-combing.   
It was mentioned in  the last section of \cite{MT2} that 
the existence of a weakly tame 1-combing actually implies the existence of a tame 1-combing. 
\end{remark}

\section{Proof of Theorem \ref{main1}}

\subsection{Comparison of qsf and wgsc conditions}
The subject of this section is the proof of  
the almost-equivalence of qsf and wgsc conditions from Theorem \ref{main1}.
Our result is slightly more general and includes the 1-tameness and 
small content conditions, which will be used later, in section 6.

\begin{proposition}
A wgsc polyhedron  has small content and is 1-tame. A polyhedron
which is either 1-tame or else has small content is qsf.
\end{proposition}
\begin{proof}
Let $C$ be a compact sub-polyhedron of the polyhedron $X$.

\vspace{0.1cm} (1). Assume that $X$ is wgsc. Then one can embed
$C$ in a compact 1-connected sub-polyhedron $K\subset X$. Taking 
then $D=E=K$ one finds that $X$ has small content and is 1-tame.

\vspace{0.1cm} (2).  Suppose that $X$ has small content, and $D$ and
$E$ are the sub-polyhedra provided by definition \ref{small}. 
Let $\gamma$ be a loop in $E$, based at a point in $C$. 
We consider the  decomposition of $\gamma$ into maximal arcs $\gamma[j]$
which are (alternatively) contained either in $D$ or in the closure 
$\overline{E-D}$ of $E-D$, namely 
$\gamma[1]\subset D, \gamma[2]\subset \overline{E-D}$, and so on. 
Thus $\gamma[2k]\subset \overline{E-D}$ has its endpoints in $D$. 
By hypothesis there exists another arc $\lambda[k]\subset D-C$
that joins the endpoints of $\gamma[2k]$. The composition  
$\gamma^k=\gamma[2k]\lambda[k]^{-1}$ is then a loop in $E-C$. 
Moreover, the composition 
$\gamma^0=
\gamma[1]\lambda[1]\gamma[3]\lambda[2]\cdots \gamma[2k-1]\lambda[k]\cdots$  
is a loop contained in $D$.  
Next $\gamma^j$ 
(based at one endpoint of $\gamma[2j]$ from $\overline{E-D}$) is homotopic
within $X-C$ to a loop $\bar{\gamma}^j\subset D-C$.

\vspace{0.1cm} \noindent Assume now  that we chose a system of
generators $\gamma_1,\dots,\gamma_n$ of $\pi_1(E)$. We will do the
construction above for each loop $\gamma_j$, obtaining the loops
$\gamma_j^k$ in
$E-C$ which are homotopic to $\bar{\gamma}_j^k$ in $D-C$. We
define first a polyhedron $\widehat{E}$ by adding to $E$ 2-disks along
the composition of the loops $\gamma_j^k(\bar{\gamma}_j^k)^{-1}$. Recall
that these two loops have the same base-point (depending on $j, k$) 
and so it makes sense to consider their composition. 

\vspace{0.1cm} \noindent There is defined a natural map
$F:\widehat{E}\to X$, which extends the inclusion
$E\hookrightarrow X$, as follows. 
There exists a homotopy  within $X-C$ keeping
fixed the base point of $\gamma_j^k$ between $\gamma_j^k$ and
$\bar{\gamma}_j^k$. Alternatively, there exists a free
null-homotopy of the loop $\gamma_j^k(\bar{\gamma}_j^k)^{-1}$ within
$X-C$. We send then the 2-disk of $\hat{E}$ capping off
$\gamma_j^k(\bar{\gamma}_j^k)^{-1}$ onto the image of the associated
free null-homotopy.

\vspace{0.1cm} \noindent It is clear that $F$ is a homeomorphism
over $C$, since the images of the extra 2-disks are disjoint from
$C$. Moreover, we claim that $\widehat{E}$ is simply connected. In
fact, any loop in $\widehat{E}$ is homotopic to a loop within $E$,
and hence to a composition of $\gamma_j$. Each $\gamma_j$ is 
homotopic rel. base point, by a homotopy in $E$, 
to $\gamma_j[1] \gamma_j^1\lambda_j[1]\gamma_j[3]\gamma_j^2\lambda_j[2]\cdots$, 
which is homotopic rel. base point,  by a homotopy in $\widehat{E}$, to 
$\gamma_j[1] \bar{\gamma}_j^1\lambda_j[1]\gamma_j[3]\bar{\gamma}_j^2\lambda_j[2]\cdots$, a loop in $D$. By hypothesis, this last loop is 
null-homotopic in $E$. Therefore
$\pi_1(\widehat{E})=0$.

\vspace{0.1cm} (3). Suppose now that $X$ is 1-tame. Let $K$ be the
compact associated to an arbitrarily given compact $C$. Any loop
$\gamma$ in $K$ is freely homotopic to a loop $\bar{\gamma}$ in
$K-C$. Consider $\gamma_1,\cdots,\gamma_n$ a system of generators
of $\pi_1(K)$. From $K$ we construct the polyhedron $\widehat{K}$ by adding
2-disks along the loops $\bar{\gamma}_j$. There exists a map
$F:\widehat{K}\to X$, which extends the inclusion
$K\hookrightarrow X$, defined as follows. The 2-disk capping off
the loop $\bar{\gamma}_j$ is sent into the null-homotopy of
$\bar{\gamma}_j$ within $X-C$. Then $F$ is obviously a
homeomorphism over $C$. Meanwhile, $\widehat{K}$ is simply
connected since we killed all homotopy classes of loops from $K$.
\end{proof}

\begin{proposition}\label{qsfe}
If the  open $n$-manifold $M^n$ is qsf and $n\geq 5$ then $M^n$ is wgsc.
\end{proposition}
\begin{proof}
It suffices to prove that any compact codimension zero 
sub-manifold $C$ is contained in a simply connected compact sub-space 
of $M^n$. By hypothesis there exists a compact connected 
and simply connected simplicial complex $K$ and a map $f:K\to M^n$ 
such that $f:f^{-1}(C)\to C$ is a PL homeomorphism.
Assume that $f$ is simplicial, after subdivision.  
Let $L$ be the 2-skeleton of $\overline{K\setminus f^{-1}(C)}$ and 
denote by $\partial L=L\cap \partial f^{-1}(C)$. Notice that 
$f^{-1}(C)$ is a manifold. 

\vspace{0.1cm} \noindent 
The restriction of $f|_L:L\to \overline{M^n\setminus C}$ to 
the sub-complex $\partial L\subset L$ is an embedding. 
Since the dimension of $L$ is 2 and $n\geq 5$, 
general position arguments show that we can perturb $f$ by a homotopy 
which is identity on $\partial L$ to a simplicial map 
$g:L\to \overline{M^n\setminus C}$ which is an embedding. 

\vspace{0.1cm} \noindent 
Observe now that $\pi_1(f^{-1}(C)\cup_{\partial L} L)\cong \pi_1(K)=0$, 
and thus $\pi_1(C\cup_{f(\partial L)} g(L))=0$. 
Take a small regular neighborhood $U$ of $C\cup_{f(\partial L)} g(L)$ 
inside $M^n$. Then $U$ is a simply connected compact sub-manifold 
of $M^n$ containing $C$. 
\end{proof}

\begin{remark}
A similar result was proved in \cite{FG} for Dehn exhaustibility. 
In particular a
$n$-manifold which is Dehn-exhaustible is wgsc provided  that $n\geq 5$.
\end{remark}

\begin{definition}
A finitely generated group has  the topological property $A$ 
if {\em some} Cayley complex has property $A$. The topological properties
$A$ and $B$ are {\em almost-equivalent} for finitely presented groups 
if a  finitely presented group has $A$ if and only if it has $B$.  
\end{definition}

\begin{corollary}
The wgsc, gsc, qsf, Dehn-exhaustibility, 1-tameness and
small content  are almost-equivalent for finitely presented groups.
\end{corollary}

\noindent This also yields the following  geometric characterization of the
qsf:

\begin{corollary}
The group $\Gamma$ is qsf if and only if the universal covering
$\widetilde M^n$ of any compact manifold $M^n$ with
$\pi_1(M^n)=\Gamma$ and dimension $n\geq 5$ is wgsc (or gsc). In
particular, a qsf group admits a presentation whose Cayley complex
is wgsc.
\end{corollary}
\begin{proof}
The ``if'' part is obvious. 
Assume then that $\Gamma$ is wgsc and thus there exists a compact 
polyhedron whose universal covering is wgsc and hence qsf. It is known 
(see \cite{BM1}) that the qsf property does not depend on the 
particular compact polyhedron we chose. Thus, if $M^n$ 
is a compact manifold with fundamental group $\Gamma$ then 
$\widetilde{M^n}$ is also qsf. By the previous Proposition, when 
$n\geq 5$ $\widetilde{M^n}$ is also wgsc, as claimed. 

\vspace{0.1cm} \noindent 
Further, if the group $\Gamma$ is qsf then consider a compact 
$n$-manifold $M^n$ with fundamental group $\Gamma$ and $n\geq 5$. 
It is known that $\widetilde{M^n}$ is qsf and thus wgsc. 

\vspace{0.1cm} \noindent 
Consider a triangulation of $M^n$ 
and $T$ a maximal tree in its 1-skeleton. Since the finite tree $T$ 
is collapsible it has a small neighborhood  $U\subset M^n$ 
homeomorphic to the $n$-dimensional disk. 
The quotient $U/T$ is homeomorphic to 
the $n$-disk and thus to $U$. This implies that the quotient $M^n/T$ is 
homeomorphic to $M^n$. Therefore we obtain a finite 
CW-complex $X^n$ homeomorphic to $M^n$ and having a single vertex. 
Also $X^n$ is wgsc since it is homeomorphic to a wgsc space.

\vspace{0.1cm} \noindent 
The wgsc property is inherited by the 2-skeleton, namely a locally finite 
CW-complex $X$ is wgsc if and only if its 2-skeleton is wgsc.
This means that the universal covering of the 2-skeleton of $X^n$ 
is wgsc. But any finite CW-complex of dimension 2 
with one vertex and fundamental 
group $\Gamma$ is the Cayley complex associated to a suitable presentation of 
$\Gamma$. Thus the Cayley complex of this presentation is wgsc, as claimed.  
 \end{proof}

\subsection{Qsf and tame 1-combability}
The subject of this section is to end the proof of Theorem \ref{main1} 
by proving that qsf and tame 1-combability are almost-equivalent for 
finitely presented groups.

\vspace{0.1cm}\noindent
We will consider below open connected manifolds with {\em finitely many 
1-handles}, which slightly generalize the gsc condition. In the smooth 
category this means that there is a proper Morse function with only  
finitely many index 1 critical points. In the PL category we  can ask 
that the manifold have 
a proper handlebody decomposition for which 1-handles and 2-handles 
are in cancelling position for all but finitely many pairs.

\begin{proposition}
Let $W^n$, $n\geq 5$, be an open connected manifold admitting a  
proper handlebody decomposition with only finitely many 1-handles. 
Then $W^n$ is Tucker. 
\end{proposition} 
\begin{proof}
We have to prove that for sufficiently large compact sub-complexes 
$K$ the group $\pi_1(W-K)$ is finitely generated. 

\vspace{0.1cm}\noindent 
Consider a proper handlebody decomposition with a single 0-handle and finitely many 1-handles. 
We shall assume that $K$ is large enough to include all index 1 handles. Further,  by compactness
there is a union of handles $C$ containing $K$. Here $C$ is a manifold with boundary 
$\partial C$. We obtain $W-{\rm int}(C)$ from $\partial C\times [0,1]$ by adding inductively 
handles of index at least  2. In particular $W-{\rm int}(C)$ has as many connected components 
as $\partial C$.   Let $F$ be a connected  component of $\partial C$ and 
$Z$ be the corresponding connected component of $W-{\rm int}(C)$.  
\begin{lemma}
The inclusion $F\hookrightarrow Z$ induces a surjective homomorphism $\pi_1(F)\to \pi_1(Z)$. 
\end{lemma}
\begin{proof}
Let $Z_k$ be the result of adding the next $k$ handles of the 
decomposition to $F\times [0,1]$ and  let $F_k$ be the other boundary of  
$Z_k$, namely $\partial Z_k=F\cup F_k$, and $Z_0=\emptyset$.  
We claim first that 
$\pi_1(Z_{k+1}-{\rm int}(Z_k), F_k)=0$, or equivalently, the 
homomorphism induced by inclusion 
$\pi_1(F_k)\to \pi_1(Z_{k+1}-{\rm int}(Z_k))$ is surjective.  
In fact one obtains $Z_{k+1}-{\rm int}(Z_k)$ from 
$F_k\times [0,1]$ by adding one handle of index at least 2.
Then Van Kampen implies the claim.  Further  
$Z_{k}=\cup_{j=0}^{k-1} (Z_{j+1}-{\rm int}(Z_j))$ so that 
$\pi_1(Z_{k})$ is the  iterated amalgamated product 
\[ \pi_1(Z_1)*_{\pi_1(F_1)}\pi_1(Z_2- {\rm int}(Z_1))*_{\pi_1(F_2)}
   \pi_1(Z_3- {\rm int}(Z_2))*\cdots *_{\pi_1(F_{k-1})} 
\pi_1(Z_k- {\rm int}(Z_{k-1}))\]
The previous claim shows then that the inclusion $F\to Z_k$ induces a 
surjective map $\pi_1(F)\to \pi_1(Z_k)$, for each $k$.  Letting $k$ go 
to infinity  we find that $\pi_1(F)$ surjects onto $\pi_1(Z)$. 
\end{proof}

\vspace{0.1cm}\noindent 
At last $W-K$ is obtained by gluing $W-C$ and $C-K$. It is clear that 
$C-K$ has finitely generated fundamental group. As $\partial C$ has 
finitely many connected components the use of Van Kampen and the lemma above 
imply that $\pi_1(W-K)$ is finitely generated. 
\end{proof}

\begin{corollary}
If $W^n$ is open gsc manifold of dimension $n\geq 5$ then $W^n$ is Tucker. 
\end{corollary}

\begin{proposition}
A finitely presented group is qsf iff it  is tame 1-combable. 
\end{proposition}
\begin{proof}
The ``if'' implication is proved in \cite{MT2}. 
Let $G$ be qsf. Choose some closed  triangulated 5-manifold $M$ with $\pi_1(M)=G$. 
According to our previous result $\widetilde M$ is an open gsc manifold.  
In particular, by the corollary above $\widetilde{M}$ is Tucker. 
Now, a complex $X$ is Tucker if and only if 
its 2-skeleton $X^{2}$ is Tucker.  
Therefore the 2-skeleton of $\widetilde{M}$ and hence 
the universal covering $\widetilde{M^{(2)}}$ of the 2-skeleton of the triangulation of $M$ 
is Tucker. It is clear that $\pi_1(M)=\pi_1(M^{(2)})$.  
Recall then from \cite{MT2} that  $G$ is tame combable if there is 
some finite 2-complex $X$ with $\pi_1(X)=G$ for which $\widetilde{X}$ has the Tucker property. 
This proves that $G$ is tame 1-combable. 
\end{proof}

\vspace{0.1cm}\noindent 
In particular we obtain the Corollary \ref{qi}, which we restate  here 
for the sake of completeness: 

\begin{corollary}
The qsf property is a quasi-isometry invariant of finitely presented groups. 
\end{corollary}
\begin{proof}
Brick proved (see \cite{Brick}, and also the refinement from 
\cite{HM2}, Theorem A) that a group quasi-isometric to a finitely 
tame 1-combable group is also tame 1-combable. 
\end{proof}

\subsection{Some examples of qsf groups}

\subsubsection{General constructions.}
It follows from \cite{Brick,BM1,BM2,MT2} that most geometric examples 
of groups are actually qsf.

\begin{example}
\begin{enumerate}
\item A  group $G$ is qsf if and only if a finite index subgroup $H$ of $G$  is qsf.
\item Let $A$ and $B$ be  finitely presented qsf groups and $C$ be a common 
finitely generated subgroup.
Then the amalgamated free product $G=A \ast _C B$ is qsf. 
If $A$ is a finitely presented qsf group and $\phi:C_1\to C_2$ is an 
isomorphism of finitely generated subgroups of $A$, then 
the HNN-extension $A \ast _\phi $ is qsf. 
Conversely, if $A,B$ are finitely presented
and $C$ is finitely generated then $A\ast _C B$ (respectively 
$A \ast _\phi$, where $\phi:C_1\to C_2$ is an 
isomorphism of finitely generated subgroups of $A$) 
is qsf implies that $A$ and $B$ are qsf. 
\item All one-relator groups are qsf.
\item The groups from the class ${\mathcal C}_+$ (combable)
in the sense of Alonso-Bridson (\cite{AB}) are qsf. In particular
automatic groups, small cancellation groups, semi-hyperbolic
groups,  groups acting properly co-compactly on Tits buildings of
Euclidean type, Coxeter groups, fundamental groups of closed
non-positively curved 3-manifolds are qsf. Notice that all these
groups have solvable word problem.
\item If a group has a tame 1-combing 
then it is qsf. In particular, asynchronously automatic groups
(see \cite{MT2}) are qsf.
\item Groups which are simply connected at infinity
are qsf (\cite{BM1}).
\item Assume that $1\to A\to G\to B\to 1$ is a short exact
sequence of infinite finitely presented groups. Then $G$ is 
qsf (\cite{BM2}). More generally, graph products 
(i.e. the free product of vertex groups with additional relations 
added in which elements of adjacent vertex groups commute with each other) 
of infinite finitely presented groups associated to nontrivial 
connected graphs are qsf.  
\end{enumerate}
\end{example}

\begin{remark}
The last property above is an algebraic analog of the fact that
the product of two open simply-connected manifolds is gsc.
Moreover, if one of them is 1-ended then the product is simply
connected at infinity.
\end{remark}

\begin{remark}
There exist finitely presented qsf groups with 
unsolvable word problem. Indeed, in \cite{CM} the authors constructed 
a group with unsolvable word problem that
 can be obtained from a free group by applying three 
successive HNN-extensions with finitely generated free 
associated subgroups.  Such a group is
 qsf from  (2) of the Example above.
\end{remark}

\subsubsection{Baumslag-Solitar groups: not simply connected at infinity}
The Baumslag-Solitar groups are given by the 1-relator presentation
\[B(m,n)=\langle a,b | ab^ma^{-1}=b^n\rangle, \; m,n\in \mathbb Z\]
Since they are 1-relator groups they are qsf. 
It is known that $B(1,n)$ are amenable, metabelian groups which
are neither lattices in 1-connected solvable real Lie groups nor
CAT(0) groups (i.e. acting freely co-compactly on a proper CAT(0)
space). 

\vspace{0.1cm} \noindent 
Notice that $B(1,n)$ are
not almost convex  with respect to any generating set and not
automatic either, if $n\neq \pm 1$. 

\vspace{0.1cm} \noindent 
Recall that a  group $G$
which is simply connected at infinity should satisfy
$H^2(G,\mathbb Z G)=0$. Since this condition is not satisfied by
$B(1,n)$, for $n > 1$ (see \cite{Mi2}),  
these groups are not simply connected at
infinity. 

\vspace{0.1cm} \noindent The higher Baumslag-Solitar groups
$B(m,n)$ for $m,n >1$ are known to be nonlinear, not residually
finite, not Hopfian (when $m$ and $n$ are coprime), 
not virtually solvable. Moreover, they are
not automatic if $m\neq \pm n$, but they are asynchronously
automatic.

\subsubsection{Solvable groups: not CAT(0)} Let $G$ be a finitely presented 
solvable group whose derived series is
\[G\triangleright G^{(1)}\triangleright G^{(2)}\triangleright\cdots
\triangleright G^{(n)}\triangleright G^{(n+1)}=1 \] If $G^{(n)}$
is finite then  $G$ is qsf if and only if 
the solvable group $G/G^{(n)}$ (whose derived length is one 
unit smaller than $G$) is qsf. 
Solvable  groups with infinite finitely generated center are qsf by 
Example 3.5.(7). More generally, if  
$G^{(n)}$ has an element of infinite order, then Mihalik 
(see \cite{Mi2}) proved that either $G$ is
simply connected at infinity or else there exist two groups
$\Lambda\triangleleft G$,  which is normal of finite index, and
$F\triangleleft \Lambda$, which is a normal finite subgroup, such that
$\Lambda/F$ is isomorphic to a Baumslag-Solitar group $B(1,m)$.
This  implies that $\Lambda$ is qsf and hence $G$ is  qsf. 
This is useful in understanding that qsf groups are far more 
general than  groups acting properly co-compactly  and by isometries 
on CAT(0) spaces. In fact, every solvable subgroup of  
such a CAT(0) group should be virtually abelian. Thus all solvable 
groups that are not virtually abelian are not CAT(0) and many of them are 
qsf (e.g. if their center is not torsion). 
Remark that there exist solvable groups with infinitely generated 
centers, as those constructed by Abels (see \cite{Abe}). 
In general we do not know whether all solvable groups (in particular those 
with finite centers) are qsf, but 
one can prove that Abels' group is qsf since it is an S-arithmetic group.

\subsubsection{Higman's group: acyclic examples} 
The first finitely presented acyclic group
was introduced by G.Higman  in \cite{Hig}:
\[ H=\langle x,y,z,w | x^w=x^2, y^x=y^2, z^y=z^2,w^z=w^2 \rangle,
\,{\rm where }\; a^b=bab^{-1}\] It is known (see e.g.\cite{DV}) 
that $H$ is an iterated amalgamated product
\[ H=H_{x,y,z}*_{F_{x,z}} H_{z,w,x},\,{\rm with }\;
H_{x,y,z}=H_{x,y}*_{F_{y}}H_{y,z}\] where $H_{x,y}=\langle x,y,
y^x=y^2\rangle$ is the Baumslag-Solitar group $B(1,2)$ in the
generators $x,y$. Here $F_{y}, F_{x,z}$ are the free groups in the
respective generators. The morphisms $F_y\to H_{x,y}$, $F_{x,z}\to
H_{x,y,z}$ and their alike are tautological i.e. they send each left
hand side generator into the generator denoted by the same letter on
the right hand side. Remark that these homomorphisms are injective. The
example above implies that $H$ is qsf. Observe that $H$ is not
simply connected at infinity  according to \cite{Mi3}. There are more general Higman groups 
$H_n$ generated by $n$ elements with $n$ relations as above in cyclic 
order. It is easy to see that $H_3$ is trivial and the arguments above 
imply that $H_n$ are qsf for any $n\geq 4$.

\subsubsection{The Gromov-Gersten examples} 
A slightly related  class of groups 
was considered by Gersten and Gromov (see \cite{Gr2}, 4.C$_3$), as follows: 
\[ \Gamma_n=\langle a_0,a_1,\ldots, a_n| a_0^{a_1}=a_0^2, \, 
a_1^{a_2}=a_1^2, \, \ldots, a_{n-1}^{a_n}=a_{n-1}^2 \rangle \]
Remark that $H_n$ is obtained from $\Gamma_n$ by adding one more 
relation that completes the cyclic order. As above $\Gamma_{n+1}$ 
is an amalgamated product $\Gamma_n*_{\Z} B(1,2)$ and thus 
$\Gamma_n$ is qsf for any $n$. These examples are very instructive since 
Gromov and Gersten proved that the connectivity radius of $\Gamma_{n+1}$ 
is an $n$-fold iterated exponential (see the next section for a discussion). 
Moreover, $\Gamma_n$ is contained in the group 
\[ \Gamma_*=\langle a_0, b| a_0^{a_0^b}=a_0^2\rangle \]
Therefore, the connectivity radius  of $\Gamma_*$ is higher 
than any iterated exponential.  
Since $\Gamma_*$ is a 1-relator group it is  qsf and has 
solvable word problem.

\subsubsection{Thompson groups} Among the first examples of infinite
finitely presented simple groups are those provided by R.Thompson in
the sixties. We refer to \cite{CFP} for a thorough introduction to
the groups usually denoted $F$, $T$ and $V$. These are by now
standard test groups.

\vspace{0.1cm} \noindent According to (\cite{BG,FH}) $F$ is  a
finitely presented group which is an ascending HNN extension of
itself. A result of Mihalik (\cite{Mi1}, Th.3.1) implies that $F$
is simply connected at infinity and thus qsf.

\begin{remark}
Notice that $F$ is a non-trivial extension of its abelianization
$\mathbb Z^2$ by its commutator $[F,F]$, which is a simple group.
However $[F,F]$ is not finitely presented, although it is still a
diagram group, but one associated to   an infinite semi-group
presentation. Thus one cannot apply directly (7) of the Example above. 
\end{remark}

\vspace{0.1cm} \noindent Moreover, the  truncated 
complex of bases due to Brown and Stein (see \cite{Br}) 
furnishes a contractible complex acted upon freely co-compactly 
by the Thompson group $V$. The 
start-point of the construction is a 
complex associated to a directed poset which is therefore 
exhausted by finite  simply connected (actually contractible) sub-complexes. 
The qsf is preserved  through all the subsequent steps of the construction
and thus the complex of bases is qsf. In particular the Thompson group 
$V$ is qsf. We skip the details.  

\begin{remark}
It is likely that all diagram groups (associated to 
a finite presentation of a finite semi-group)  considered by Guba and Sapir 
in \cite{GS1} (and their generalizations, the picture groups) are qsf. 
Farley constructed in \cite{Fa1}  free proper actions by 
isometries of diagram groups on infinite dimensional 
CAT(0) cubical complexes. This action is not co-compact and moreover the
respective CAT(0) space is infinite dimensional. However there
exists a natural construction of truncating the  CAT(0) space $X$  in order
to get subspaces $X_n$ of $X$ which are invariant, co-compact
and $n$-connected. Farley's construction works well (\cite{Fa2})
for circular and picture diagrams (in which planar diagrams are
replaced by annular diagrams or diagrams whose wires are crossing
each other). However, these sub-complexes are not 
anymore CAT(0), and it is not clear whether they are qsf. 
Notice that these groups have solvable word problem (see \cite{GS1}). 
\end{remark}
  




\subsubsection{Outer automorphism groups.}
In the case of surface groups these correspond to mapping class groups. 
Since a finite index subgroup acts freely 
properly discontinuously on the Teichm\"uller space it follows that 
mapping class groups are qsf. The study of Morse type functions on the 
outer space led to the fact that $Out(F_n)$ is $2n-5$ connected 
at infinity and thus qsf as soon as $n\geq 3$ (see \cite{BF}).

\section{Proof of Theorem \ref{main2} and applications}

\subsection{The qsf growth}
Let $P$ be a  finite presentation of the qsf group  $G$ and
$C(G,P)$ be the associated Cayley complex. 
There is a natural word metric on the set of vertices of the 
Cayley graph $C^1(G,P)$ (the 1-skeleton of the Cayley complex) by setting 
\[ d(x,y) =\min |w(xy^{-1})|\]
where $|w(a)|$ denotes the length of a word $w(a)$ in the letters 
$s, s^{-1}$, for $s\in P$, representing the element $a$ in the group $G$. 
By language abuse we call metric complex a simplicial complex whose 
0-skeleton is endowed with a metric. 

\begin{definition}
The metric ball $B(r,p)\subset C(G,P)$ (respectively metric sphere 
$S(r,p)\subset  C(G,P)$) of radius $r\in\Z_+$ centered at 
some vertex $p$ is the following sub-complex of $C(G,P)$:
\begin{enumerate}
\item the vertices of $B(r,p)$ (respectively $S(r,p)$) 
are those vertices of $C(G,P)$ 
staying at distance at most $r$ (respectively $r$) from $p$;
\item the edges and the 2-cells of $B(r,p)$ (respectively $S(r,p)$) 
are those edges and 2-cells 
of $C(G,P)$  whose boundary vertices are at distance at most $r$ 
(respectively $r$) from $p$. 
\end{enumerate}  
Denote by $B(r)$ (respectively $S(r)$) the metric ball (respectively  sphere) 
of radius $r$ centered at the identity. 
\end{definition}

\begin{definition}
A {\em $\pi_1$-resolution} of the
polyhedron $C$ inside $X$ is a pair $(A,f)$, where $A$ is a CW complex 
and $f:A\to X$ a PL map such that
$f:f^{-1}(C)\to C\subset X$ is a PL-homeomorphism and $\pi_1(A)=0$. 
\end{definition}

\vspace{0.1cm}
\noindent 
We want to refine the qsf property  for metric complexes. 
As we are interested in 
Cayley complexes below we formulate the definition in this context: 
\begin{definition}
The {\em qsf growth} function $f_{G,P}$ of the Cayley complex $C(G,P)$, 
is:
\[ f_G(r)=\inf \{ R  {\rm \; such\; that\; there\; exists\; a}\; \pi_1{\rm -resolution\;
of }\; B(r)\; {\rm into}\; B(R)\}\] 
\end{definition}

\vspace{0.1cm}
\noindent 
Recall that  the real functions $f$ and $g$ are {\em rough equivalent} 
if there exist constants $c_i, C_j$ (with
$c_1,c_2 >0$) such that
\[ c_1 f(c_2 R)+ c_3 \leq g(R) \leq C_1 f(C_2 R) + C_3 \]
One can show easily that the
rough equivalence class of $f_{G,P}(r)$ depends only on the group $G$
and not on the particular presentation, following \cite{BM1} and \cite{FO}.
We will write it as $f_G(r)$.  
We don't know whether the rough equivalence 
class of $f_G$ is a quasi-isometry invariant. 
This would be true if we could compare 
$f_G$ with the tameness function 
of Hermiller and Meier (\cite{HM2}).

\vspace{0.1cm} \noindent Recall from (\cite{Gr2}, 4.C) that the 
{\em connectivity radius} $R_1(r)$  defined by Gromov 
is the infimal $R_1(r)$ such that $\pi_1(B(r))\to \pi_1(B(R_1(r))$ 
is zero. Notice that the rough equivalence class of $R_1$ 
is  also well-defined and independent on the group presentation we chose 
for the group.

\begin{remark}
Observe that $\pi_1(B(r))\to \pi_1(B(f_G(r)))$ is zero. Thus $f_G$ is
bounded from below by the connectivity radius $R_1$. 
\end{remark}

\noindent Recall that the {\em isodiametric function} of a group $G$, following
Gersten, is the infimal $I_G(k)$ so that loops of length $k$ bound
disks of diameter at most $I_G(k)$ in the Cayley complex. The rough
equivalence class of $I_G$ is a quasi-isometry invariant of the
finitely presented group $G$.

\begin{proposition}
A qsf group whose qsf growth  $f_G$ is recursive
has a solvable word problem.
\end{proposition}

\begin{proof}
Observe that the growth rate of the qsf is
an upper bound for the Gersten isodiametric function $I_G$, and 
the word problem is solvable whenever the isodiametric function 
is recursive. This is standard: if a word $w$ is trivial in the group 
$G$ presented as $G=\langle S|R\rangle$, 
then it is a product of conjugates of relators 
$uru^{-1}$, $r\in R$. By the definition of the isodiametric function 
one can choose these conjugates in such way that 
$|u|\leq I_G(|w|) \leq f_G(|w|)$ and 
this leads to a finite algorithm that checks whether $w$ is trivial or not. 
\end{proof}

\subsection{Metric balls and spheres in Cayley complexes}
We consider now some  metric complexes  satisfying a closely related 
property. On one side this condition seems to be slightly weaker than the wgsc 
since we could have nontrivial (but uniformly small) loops, but on the other 
side the exhaustions we consider are restricted to metric  balls. 

\begin{definition}
A metric complex has {\em $\pi_1$-bounded balls} (respectively {\em spheres}) if
there exists a constant $C$ so that $\pi_1(B(r))$ (respectively 
$\pi_1(S(r))$) is normally generated by loops with length smaller
than $C$.
\end{definition}
\begin{remark}
If the balls in a metric complex are simply connected then 
the complex is obviously wgsc. However, if the complex is wgsc it is not 
clear whether we can choose an exhaustion by simply connected metric balls.
Thus the main constraint in the definition above is the requirement 
to work with metric balls.  
\end{remark}
\noindent We actually show that this condition   
puts strong restrictions on the
group:
\begin{proposition}
If some Cayley complex of a finitely presented group has $\pi_1$-bounded 
balls (or spheres) then the group is qsf with  linear qsf growth.
\end{proposition}
\begin{proof}
Any loop $l$ in the ball $B(r)$ is null-homotopic in the Cayley complex. Thus
there exists a disk-with-holes  lying in $B(r)$ such that 
the outer boundary component is equal to $l$ and the other boundary 
components $l_1',\ldots l_p'$ lie in $S(r)$. 

\vspace{0.1cm} \noindent
If we have $\pi_1$-bounded  spheres
then we can assume that each loop $l_j'$, $1\leq j\leq p$,  is  
made of uniformly small loops on $S(r)$ connected by means of arcs. 
A loop of length $C$
in the Cayley graph bounds a disk of diameter $I_G(C)$ in the 
Cayley complex. Thus these disks 
have uniformly bounded diameters. The loop $l_j'$ bounds therefore a
disk $D(l_j')$ which is disjoint from $B(r-I_G(C))$ and lies within
$B(r+I_G(C))$.

\vspace{0.1cm} \noindent We can use this procedure for a 
system of loops $l_j$, $1\leq j\leq n$ which
generate $\pi_1(B(r))$. Thus, to any loop $l_j$ we associate a disk-with-holes  
having one boundary component $l_j$  while the other boundary 
components are the loops $l_{j,k}'$ which lie on $S(r)$ 
and then null-homotopy disks $D(l_{j,k}')$ 
as above. Let $D_j$ denote their union, which is a 2-dimensional sub-complex 
of $C(G,P)$ providing a null-homotopy of $l_j$. 

\vspace{0.1cm} \noindent 
Let $A$ denote the simplicial complex made of $B(r)$ 
union a number of 2-disks $D(j)$ which are attached  to $B(r)$ 
along the loops $l_j$. As the set of loops $l_j$ generate $\pi_1(B(r))$ 
the complex $A$ is simply connected. 

\vspace{0.1cm} \noindent
We define the map $A\to
B(r+I_G(C))$ by sending each disk $D(j)$ into the corresponding  null-homotopy 
disk $D_j$. Since $\pi_1(A)=0$ this  map provides  a 
$\pi_1$-resolution of $B(r-I_G(C))$.

\vspace{0.1cm} \noindent The same proof works for $\pi_1$-bounded 
balls.
\end{proof}

\vspace{0.1cm} \noindent 
{\em End of the proof of Theorem \ref{main2}.} We have to show that 
a group having a Cayley complex with  $\pi_1$-bounded balls or spheres
has linear connectivity radius and solvable word problem.
The connectivity radius is at most linear since loops generating 
$\pi_1(B(r))$ are null-homotopic using uniformly bounded null-homotopies 
whose size depends only on $C$. 
Thus $\pi_1(B(r))\to \pi_1(B(r+C'))$ is zero for $C'\geq I_G(C)$.  
This means that a loop in
$B(r)$ is null-homotopic in the Cayley complex only if it is null-homotopic 
within $B(r+C')$. Moreover, the last condition can be checked by a 
finite algorithm for given $r$, and in particular one can check 
whether a given word of length $r$ is trivial or not.

\begin{remark}\label{smallb}
Some Cayley complexes of hyperbolic groups have $\pi_1$-bounded balls and
spheres. For instance this is so for any of the Rips 
complexes, whose metric balls are known to be simply connected. 
It is likely that any Cayley complex associated to a finite presentation 
of a hyperbolic group has $\pi_1$-bounded balls. 
Furthermore, if a group  $G$ acts properly co-compactly on  
a CAT(0) space then the metric balls are convex  and thus they are 
simply connected. It seems that this implies that any other space that is 
acted upon by the group $G$ properly co-compactly  (thus quasi-isometric to 
the CAT(0) space) should have also $\pi_1$-bounded balls.  
This would follow if the $\pi_1$-bounded balls property were a quasi-isometry 
invariant. 
\end{remark}

\begin{remark}
One can weaken the requirements in  the definition of
$\pi_1$-bounded spheres, in the case of a Cayley complex of a group,
as follows. We only ask that the group 
$\pi_1(B(r))$ be normally generated by loops
of length $\rho(r)$ where
\[ \lim_{r\to\infty}r-I_G(\rho(r)) =\infty \]
Note that the limit should be infinite for any choice of the 
isodiametric function $I_G$
within its rough equivalence class.  Then, finitely presented groups 
verifying this weaker condition are also qsf, by means of the same 
proof.

\vspace{0.1cm} \noindent Notice however that
$I_G(r)$ should be non-recursive for groups  $G$ with non-solvable word
problem, so that $\rho(r)$  grows extremely slow if non-constant.
Moreover, if we  only ask that the function $r-I_G(\rho(r))$ be
recursive then the group under consideration should have again
solvable word problem. In fact, we have by the arguments above the
inequalities
\[I_G(r-I_G(\rho(r))\leq f_G(r-I_G(\rho(r)) \leq r+I_G(\rho(r)) < 2r \]
and thus $I_G(r)$ is recursive since it is bounded by the inverse of a
recursive function.
\end{remark}

\begin{remark}
Recall that the Gersten-Gromov groups $\Gamma_n$ have $n$-fold iterated 
exponential connectivity radius, and thus at least that large 
qsf growth, while $\Gamma_*$ has connectivity radius 
higher than any iterated exponential (see \cite{Gr2}, 4.C$_3$). 
We saw above that all these groups are qsf. 
However the last corollary shows that the metric balls in their 
Cayley complexes are not $\pi_1$-bounded, and thus their exhaustions by simply 
connected sub-complexes should be somewhat exotic. On the other hand 
we can infer from Remark \ref{smallb} that 
their Cayley complexes have not (group invariant) 
CAT(0)-metrics although they are both aspherical and qsf.   
\end{remark}

\subsection{Rewriting systems}
Groups admitting a rewriting system form a particular class among
groups with solvable word problem (see \cite{HM1} for an extensive
discussion). A {\em rewriting system} consists of several replacement
rules
\[ w^+_j\to w^-_j \]
between words in the generators of the presentation  
$P$. We suppose that both $s$ and $s^{-1}$ belong to $P$.  
A reduction of the word $w$
consists of a replacement of some sub-word of $w$ according to one
of the replacement rules above. The word is said irreducible if no
reduction could be applied anymore. The rewriting system  is
{\em complete} if for any word in the generators the reduction process
terminates in finitely many steps and is said to be {\em confluent} if 
the irreducible words obtained at the end of the reduction are 
uniquely defined by the class of the initial word, as an element of
the group. Thus the irreducible elements are the normal forms for
the group elements. If the rules are not length increasing then
one calls it a {\em geodesic rewriting system}. We will suppose  that
the rewriting system  consists of finitely many rules.

\begin{proposition}
A finitely presented group admitting a complete confluent geodesic rewriting
system is qsf.
\end{proposition}
\begin{proof}
In \cite{HM1} is proved that such a group is almost convex and
thus  qsf (by Proposition 4.6.3, see also \cite{O}).

\vspace{0.1cm} \noindent Here is a shorter  direct proof. We prove
that actually the balls $B(r)$ in the Cayley complex are simply
connected. Observe first that in any Cayley complex we have:
\begin{lemma}
The fundamental group $\pi_1(B(r))$ is generated by loops of
length at most $2r+1$.
\end{lemma}
\begin{proof}
Consider a loop $ep_1p_2...p_ke$ based at the identity element
$e$ and sitting in $B(r)$. Here $p_j$ are the consecutive vertices of
the loop. There exists a geodesic $\gamma_j$ that joins $p_j$ to
$e$, of length at most $r$. 
It follows that the initial loop is the product of loops
$\gamma_j^{-1}p_jp_{j+1}\gamma_{j+1}$. Since $p_j\in B(r)$ all these
loops have length at most $2r+1$.
\end{proof}
\noindent Consider now the Cayley complex of a group presentation
that includes all rules from the rewriting system. This means that
there is a relation associated to each rule $w^+_j\to w^-_j$. We
claim that the balls $B(r)$ are simply connected. By the previous
lemma it suffices to prove that loops of length at most $2r+1$
within $B(r)$ are null-homotopic in $B(r)$.

\vspace{0.1cm} \noindent  
Choose such  a loop in $B(r)$ which is
represented by the word $w$ in the generators.  We can assume that 
the normal form of the identity element is the trivial word. 
Since the loop is
null-homotopic in the Cayley complex the word $w$ should reduce to
identity by the  confluent rewriting system. Let then consider some reduction
sequence:
\[ w\to w_1\to w_2\to \cdots w_N\to e \]
Each word $w_i$ represents a loop based at the identity 
in the Cayley graph. Each step
$w_j\to w_{j+1}$ is geometrically realized as a homotopy 
in which the loop associated to the word $w_j$ is slided
across a 2-cell associated to a relation from the rewriting
system. Further the lengths of these loops verify
 $|w_j|\geq |w_{j+1}|$
since the length of each reduction is non-increasing, by
assumption. Thus $|w_j|\leq 2r+1$ and this implies that the loop
is contained within $B(r)$. This proves that the reduction
sequence above is a null-homotopy of the loop $w$ within $B(r)$.
\end{proof}

\begin{remark}
The Baumslag-Solitar groups $B(1,n)$ and the solvgroups (i.e.
lattices in the group SOL) admit rewriting system but not geodesic
ones (\cite{HM1}), since they are not almost convex.
\end{remark}
\begin{remark}
More generally one proved in \cite{HM1} that groups admitting a
rewriting system are tame 1-combable  and thus qsf
by \cite{MT2}.
\end{remark}

\begin{remark}
One might wonder whether finitely presented groups that 
have solvable word problem are actually qsf.
Notice that an algorithm solving the word problem  does not 
yield a specific null-homotopy  disk for a given loop in the 
Cayley complex, but rather checks whether a given path closes up. 
\end{remark}
\begin{remark}
The geometry of null-homotopy disks (size, diameter, area) 
is controlled by the various filling functions associated to the group. 
However, in the qsf  problem one wants to understand the 
position of the null-homotopy disks with respect to exhaustion 
subsets, which is of topological nature. 
The choice of the exhaustion is implicit but very important and it 
should depend on the group under consideration.  
\end{remark}

\section{Extensions by finitely generated groups and the Grigorchuk group}

\subsection{Infinitely presented groups.} Although it does not  
make sense to speak, in general, of the qsf property for an 
infinitely presented group, one can do it if, additionally, we specify a 
group presentation.

\vspace{0.1cm} \noindent 
Recall first that the elementary Tietze transformations of group presentations 
are the following: 

\begin{enumerate}
\item[(T1)] {\em Introducing a new generator}. One replaces 
$\langle x_1,x_2,\ldots| r_1,r_2,\ldots\rangle$ by 
$\langle y, x_1,x_2,\ldots| ys^{-1},r_1,r_2,\ldots\rangle$, where 
$s=s(x_1,x_2,\ldots)$ is an arbitrary word in the generators $x_i$. 
\item[(T2)]  {\em Canceling a generator}. This is the inverse of (T1).
\item[(T3)]  {\em Introducing a new relation}. One replaces  
$\langle x_1,x_2,\ldots| r_1,r_2,\ldots\rangle$  by 
$\langle x_1,x_2,\ldots| r, r_1,r_2,\ldots\rangle$, where 
$r=r(r_1,r_2,\ldots)$ is an arbitrary word in the conjugates of 
relators $r_i$ and their inverses.
\item[(T4)]  {\em Canceling a relation}. This is the inverse of (T3). 
\end{enumerate}

\begin{definition}
We say that two infinite presentations are {\em finitely equivalent} 
or, they belong to the same finite equivalence class, if there 
exists a finite sequence of elementary 
Tietze moves that changes one presentation 
into the other. 
\end{definition}

\begin{proposition}\label{finequiv}
The qsf property is well-defined for groups with a presentation
from a fixed finite equivalence class: 
if the Cayley complex $C(G,P)$ is qsf then the Cayley complex 
$C(G,Q)$ is qsf for any presentation $Q$ of $G$ which is finitely equivalent 
to $P$.  
\end{proposition}
\begin{proof}
The proof from \cite{BM1} works in this case word-by-word. 
\end{proof}
\noindent 
Most of the properties shared by the qsf finitely presented groups 
hold, more generally, for the qsf infinitely presented groups equipped  
with the convenient finite equivalence class of presentations. 
For instance, for any infinite groups $G, H$ and group presentations 
$(G,P_G)$ and $(H,P_H)$, 
the Cayley complex  $C(G\times H, P_G\times P_H)$ is 
qsf, where $P_G\times P_H$ is the product presentation of $G\times H$. 
     
\begin{remark}
One can obtain an infinite presentation of a group  
whose Cayley complex is not wgsc, 
by the method from Remark 2.5. However, it is more  
difficult to prove that a {\em specific} infinite presentation of some  
group is qsf, for instance in the case of Burnside groups.  
\end{remark}
\begin{remark}
The previous proposition might be extended farther. In fact one could allow 
infinitely many Tietze moves, if they do not accumulate at finite distance
but the complete definition is quite involved. 
\end{remark}

\subsection{Extensions by infinitely presented groups.} One method 
for constructing finitely presented groups is to use suitable 
extensions of finitely presented groups by infinitely presented ones. 
We did not succeed in proving that {\em all} such extensions 
are qsf. However, for finitely presented  extensions by {\em 
finitely generated} groups things  might simplify 
considerably. We start with the following definition from \cite{Bar}: 

\begin{definition}
An {\em endomorphic} presentation is an expression of the form 
$P= \langle S | Q| \Phi| R\rangle$, where 
$S$ is an alphabet (i.e. a set of symbols), $Q,R$ are sets of reduced words in the 
free group $F(S)$ generated by $S$ and $\Phi $ is a set of 
injective free group homomorphisms $F(S)\to F(S)$. The
endomorphic presentation is {\em finite}  if  
all sets $S,Q,\Phi, R$ are finite. This data defines the group: 
\[ G(P)=F(S)/\langle \, Q\cup \bigcup_{\phi\in \Phi^*}\phi(R)\,\,\rangle^{\sharp}\]
where $\langle, \rangle^{\sharp}$ denotes the normal closure and 
$\Phi^*$ is the monoid generated by $\Phi$ i.e. the closure 
of $\Phi\cup \{1\}$ under the composition. The endomorphic 
presentation is said to be {\em ascending} if $Q=\emptyset$. 
\end{definition}

\vspace{0.1cm} \noindent 
Bartholdi observed that groups with  finite 
ascending endomorphic presentations are naturally contained in 
finitely presented groups constructed as generalized ascending HNN extensions, 
by adding finitely many stable letters. 
Each  $\phi\in \Phi$ induces a group endomorphism 
$\varphi:G\to G$ and we suppose that the correspondence $\phi\to \varphi$ is 
one-to-one so $\Phi$ is also a set of endomorphisms of $G$. 
Then the finitely presented group 
\[ \overline{G}=\langle S\cup \Phi |Q\cup R \cup \{\phi^{-1} s \phi=\phi(s); s\in S, 
\phi\in \Phi\} \rangle \]
is a {\em generalized HNN extension} of $G$ with 
stable letters corresponding to $\varphi\in \Phi$. 
If the endomorphic presentation is ascending (i.e. $Q=\emptyset$), and 
the endomorphisms $\varphi:G\to G$ are injective then 
the natural homomorphism $G\to \overline{G}$ is an embedding and 
$\overline{G}$ will be what is standardly called an ascending 
HNN extension with set of stable letters $\Phi$. 
Further, if 
the natural map $G\to \overline{G}$ is an embedding, then we can assume that 
the endomorphic presentation of $G$ is ascending. In fact the relations from  
$Q$ and the conjugacy relations in $\overline{G}$ imply that 
the relations $\cup_{\phi\in\Phi^*}\phi(Q)$ are satisfied in $G$. 
Thus we can replace $R$ by $R\cup Q$ in the endomorphic presentation 
of $G$ and obtain the same generalized HNN extension group 
$\overline{G}$.

\vspace{0.1cm} \noindent
Set $N(G)$ for the normal subgroup of $\overline{G}$  
generated by $G$. We have  then an exact sequence  
\[ 1\to N(G) \to \overline{G}\to L \to 1 \]
where  the quotient $L=\overline{G}/N(G)$ has the presentation $P_L$ below: 
\[ L=\langle S\cup \Phi |S\cup R \cup \{\phi^{-1} s \phi=\phi(s); s\in S, 
\phi\in \Phi\} \rangle \] 
Using elementary Tietze moves one sees that $L$ is the free group generated by 
the  set of stable letters $\Phi$.


\vspace{0.1cm} \noindent
The images of elements of $\Phi^*\setminus\{1\}$ in $L$ will be called 
the {\em positive} elements of $L$. 
Let then $\Phi=\{\varphi_1,\varphi_2,\ldots,\varphi_k\}$.

\vspace{0.1cm} \noindent
Very interesting examples of  groups with finite ascending endomorphic 
presentations which are also branch groups appeared in the constructions 
of Bartholdi (see \cite{Bar}).

\begin{theorem}\label{ext}
Let $G$ be a finitely generated group admitting a finite ascending 
endomorphic presentation  $P_G$ such that each $\varphi_j\in \Phi$ 
is injective and $\overline{G}$ be its associated HNN extension. 
Assume that the group $G$ endowed with the presentation $P_G$ is 1-tame. 
Then $\overline{G}$ is qsf. 
\end{theorem}

\begin{remark}
The words from $\phi(R)$, $\phi\in \Phi^*$ are {\em unreduced} words 
in the free group $F(S)$, namely one can have adjacent 
canceling letters.  
This will be essential for the proof of Lemma \ref{sat}. 
Allowing unreduced words makes the hypothesis that the presentation 
$P_G$ of $G$ is 1-tame more difficult to check  and potentially more 
restrictive than in the case when the words from $R_{\infty}$ are 
reduced.  
\end{remark}

\subsection{Plan of the proof of Theorem \ref{ext}.}
Assume that we have 
a HNN extension as in the hypothesis which induces the exact sequence  
\[ 1\to N(G) \to \overline{G} \to L \to 1 \]
Consider then an infinite endomorphic presentation 
$P_G:  G=\langle A |R_{\infty} \rangle=
\langle a_i|R=\{R_j\}|\Phi=\{\phi_1,\ldots,\phi_k\} \rangle $, with an infinite 
set of relators  $R_{\infty}=\cup_{j, \phi\in \Phi^*} \phi(R_j)$ 
for $G$ and the standard presentation 
$P_L: L=\langle B | \emptyset \rangle$ for the free group $L$
with generators set $B=\{b_1,...,b_k\}$, where every $b_j$ 
corresponds to the stable letter $\phi_j$.  
One obtains an infinite presentation  for 
$\overline{G}$ by putting together the 
two presentations above, namely: 
$P_{\overline{G}}(\infty): \overline{G}=
\langle A \cup B | R_{\infty}  \cup T\rangle$, where 
the elements of $T$ express the HNN conditions for stable letters.  Thus,  
each element of $T$ has the form 
$b^{-1}ab w(a,b)$, where $a\in A$, $b\in B$ and 
$w(a,b) \in F(A)$ is some word in the generators $A$.
We call them {\em conjugacy relations}.  
Henceforth we suppose that $T$ is given by  
$T=\{ b_j^{-1} a_i b_j=\phi_j(a_i)\in A^*, a_i\in A, b_k \in B \}$.  
Notice that $B$ contains only positive letters. 
Thus it might not make sense to consider relations 
of type $b_j^{-1}a b_j$ unless $a$ is a word representing 
an element of the image of $G$ by the endomorphism associated to $b_j$.

\vspace{0.1cm}
\noindent 
However, the group $\overline{G}$ could be defined by the same set 
of generators $A\cup B$ and a finite subset of relations from above. 
We can assume that this finite presentation is 
$P_{\overline{G}}: \overline{G}=
\langle A \cup B | R \cup T\rangle$, where 
$R\subset R_{\infty}$ is a finite set of relations. 

\vspace{0.1cm}
\noindent
The plan of the proof is as follows.
The exact sequence induces a kind of foliation of the Cayley complex 
$C(\overline{G}, P_{\overline{G}})$ by horizontal leaves associated to 
$N(G)$. These leaves are connected by means of vertical tubes associated 
to conjugacy relations and going upward. Given a compact 
$C$ we can use these vertical tubes to push up loops in $C$ 
and find a larger compact $K$ whose fundamental group is generated 
by loops lying in a  top horizontal leaf $K_u$ far from $C$. 
If $G$ were finitely presented  
then the horizontal leaves would be simply connected so that loops 
could be homotopically killed inside the top leaf.  
When $G$ is not finitely presented the fundamental group of a 
(connected component of a) horizontal leaf  is generated by the relations 
in $R_{\infty}-R$. Thus loops in the top horizontal leaf are now 
freely homotopic to relation loops expressing words from  $R_{\infty}-R$. 
The 1-tameness of $C(G,P_G)$ enables us to consider only relation loops 
which are either contained in a larger compact $E$ of the horizontal 
leaf and which are disjoint from $K_u$ or else contained in $K_u$.   
Each such loop has a particularly nice 
null-homotopy in $C(\overline{G}, P_{\overline{G}})$, by expressing 
the relation as an element of $\Phi^*(R)$. Namely, there exists a canonical 
vertical tube going downward  from that loop to a 
loop which is null-homotopic in the bottom horizontal leaf. 
We add then more material to $K$ so that all canonical homotopies 
of loops from the top leaf be  either contained in $K$ or else 
disjoint from $C$. Here one makes use of the fact that the monoid of positive 
elements in $L$ defines an order on the set of horizontal leaves.  
Then  loops in $K$ are freely homotopic 
to loops which are null-homotopic in $C(\overline{G}, P_{\overline{G}})-C$ and 
so the Cayley complex is 1-tame.

\subsection{Preliminaries from Brick and Mihalik and the proof of Theorem \ref{main4}.}
Our aim is to prove that the Cayley complex 
$C(\overline{G},P_{\overline{G}})$ is qsf. 
We follow closely  the proof given by Brick and Mihalik  
in  \cite{BM2} for the fact that the extension of an infinite 
finitely presented group by an infinite finitely presented group is qsf.  
First, we state below the necessary adjustments for the main lemmas from  
\cite{BM2} work now for HNN extensions.  
Then we will point out the arguments which have to be 
modified in the present setting.  

\vspace{0.1cm}
\noindent 
Let $C(H,P_H)$ denote the Cayley 2-complex associated to  the presentation 
$P_H$ of the group $H$.
Consider now the sub-complex $X(G) \subset  C(\overline{G},P_{\overline{G}})$ 
spanned by the vertices of $N(G)\subset \overline{G}$. 
Observe that the sub-complex $X_0(G)$ spanned by $G\subset N(G)$ can be 
obtained from $C(G,P_G)$ by removing the 
2-cells corresponding to the relations from $R_{\infty}-R$. 
Moreover, $X(G)$ is the disjoint sum of copies of $X_0(G)$. 
In fact for each coset $w\in N(G)/G$ we have a copy  
$w X_0(G)\subset X(G)$ obtained by left translating by $w$. 
These copies are disjoint because any edge of $X(G)$ 
corresponds to a generator of $G$ and so a path in $X(G)$ 
corresponds to an element of $G$. Thus $wX_0(G)$ intersects $w'X_0(G)$ 
only if $wG=w'G$, for $w,w'\in N(G)$, and in this case they coincide.

\vspace{0.1cm}
\noindent 
To each $x\in L$ we associate the {\em horizontal slice} 
$x X(G) \subset C(\overline{G}, P_{\overline{G}})$ obtained by a left 
translation of $X(G)$  so that it projects down onto 
$x\in L$ under the map $\overline{G}\to L$. 
The  $0$-skeleton of the Cayley complex $C(\overline{G}, P_{\overline{G}})$  
is then decomposed as the disjoint union of  $0$-skeleta of 
slices $x X(G)$, over $x\in L$.  
The paths (edges) which are contained in such a horizontal slice 
$x X(G)$ will be called {\em $A$-paths} (respectively $A$-edges). 
The {\em $B$-edges} are those
edges of $C(\overline{G}, P_{\overline{G}})$ 
which project onto the generators $B$ of $L$. 
The 2-cells corresponding to relators in $T$ will be 
called {\em conjugation cells}. Notice that the attaching map 
of a conjugation cell is of the form $b_1ab_2^{-1}w$, where 
$b_1$ and $b_2$ are  $B$-edges corresponding to the same $b\in B$, and 
$w$ is the $A$-path corresponding to the word $w=w(a,b)$ appearing 
in the respective conjugacy relation. 
The loops in $ C(\overline{G}, P_{\overline{G}}) $ are called of {\em type 
1} if they are conjugate 
to $A$-loops and of {\em type 2} otherwise. 

\vspace{0.1cm}
\noindent 
Those sub-complexes of 
$C(\overline{G}, P_{\overline{G}})$ 
which are finite, connected and intersect each 
$w X_0(G)$, for $w\in \overline{G}$,   
in a connected --  possibly empty -- subset are called {\em admissible}.

\vspace{0.1cm}
\noindent 
Let $C\subset C(\overline{G}, P_{\overline{G}})$ be a finite connected sub-complex. 
By adding finitely many edges we may assume that $C$ is admissible. 
We want to show that there is a larger sub-complex $K$, obtained by 
adjoining finitely many edges and conjugation cells such that 
$\pi_1(K)$ is normally generated by finitely many loops in $K$
which are null-homotopic in $C(\overline{G}, P_{\overline{G}})-C$.

\begin{lemma}\label{add}
Let $Z$ be an admissible sub-complex of $ C(\overline{G}, P_{\overline{G}})$. Set 
$\{u_1,u_2,\dots,u_n\}$ for a system of generators of $\pi_1(Z)$. 
\begin{enumerate}
\item If $e$ is an $A$-edge that meets $Z$ then $Z\cup e$ is admissible. 
Further, $\pi_1(Z\cup e)$ is generated by $\{u_1,u_2,\dots,u_n,$
$\tau \lambda\tau^{-1}\}$, where $\lambda $ is an A-loop in the 
copy of $X(G)$ containing $e$. 
\item If $\Delta$ is a conjugation cell with 
$\partial \Delta=b_1^{-1}ab_2w^{-1}$, and $b_1\subset Z$ then 
$Z\cup \Delta$ is admissible. Further 
$\pi_1(Z\cup \Delta)$ is generated by 
$\{u_1,u_2,\dots,u_n,v_1,\dots,v_m\}$ where each $v_i$ is a 
type 1 generator. Moreover, if $a\subset Z$ then 
each $v_i$ is of the form $\tau \lambda\tau^{-1}$ where 
$\lambda$ is an $A$-loop in the copy of $X(G)$ containing 
the endpoint of $b_1$. 
\end{enumerate}
\end{lemma}
\begin{proof}
Lemma 2.1 from \cite{BM2} extends trivially to the present situation. 
\end{proof}

\begin{proposition}\label{2gen}
There exists a complex $C_1$ obtained from $C$ 
by adjoining finitely many $A$-edges and  conjugation cells
as in Lemma \ref{add} such that $\pi_1(C_1)$ is generated by classes of loops 
of type 1. 
\end{proposition}
\begin{proof}
We want to transform each loop of type 2 into a loop of type 1 by   
using homotopies which can be realized after 
adjoining finitely many $A$-edges and  conjugation cells to $C$ 
(satisfying the requirements of Lemma \ref{add}).
Since one does not create any additional type 2 loop we end up 
with a complex $C_1$ whose fundamental group is generated by classes 
of loops of type 1.

\vspace{0.1cm}
\noindent 
Consider first the case when there is only one stable letter, 
$B=\{b\}$. We use the conjugacy relation 
$b^{-1}a=\phi(a)b$, for $a\in A$, to move the $B$-edge labeled 
$b^{-1}$ to the right of the adjacent $A$-edge labeled $\phi(a)$. 
In meantime use the conjugacy relation  
$ab=b\phi(a)$ to move the $B$-edge  
labeled $b$ to the left of the adjacent $A$-edge labeled $\phi(a)$. 
Keep moving $B$-edges this way until two 
$B$-edges labeled $b^{-1}$ and $b$ become adjacent,  in which case 
the two edges will be removed as their labels cancel and resume the process.
This procedure eventually stops when the initial loop is transformed into 
the composition of an $A$-path with a $B$-path (or vice-versa). 
Now, the extension $\overline{G}\to L$ splits (because $L$ is free) 
and hence the $B$-path above should be a loop.   
This  $B$-loop is then homotopically trivial in its own image, 
since $L$ is free and so we eventually obtain an $A$-loop.

\vspace{0.1cm}
\noindent
Suppose now that the number of stable letters is $k\geq 2$. 
Set $\widetilde{B}=B\setminus\{b_{k}\}$, 
$\widetilde{\Phi}=\Phi\setminus\{\phi_{k}\}$, 
and $\widetilde{G}$ for the HNN extension associated to  the  set 
of stable letters $\widetilde{\Phi}$ (which is 
not necessarily finitely presented).  
Then $G$ has a natural injective homomorphism $i:G\to\widetilde{G}$ 
and thus $\varphi_{k}$ induces an isomorphism 
$\widetilde{\varphi}_{k}:i(G)\to i(\varphi_{k}(G))$ between two subgroups 
of $\widetilde{G}$. 
Therefore the group $\overline{G}$ is also the  
HNN extension  $\widetilde{G}*_{\widetilde{\varphi_{k}}}$ with base group 
$\widetilde{G}$, stable letter $b_{k}$  
and associated subgroups $i(G)$ and $i(\varphi_{k}(G))$. 

\vspace{0.1cm}
\noindent
Any loop in $C(\overline{G}, P_{\overline{G}})$ 
corresponds to a word $W$ in $\overline{G}$ representing the trivial 
element in the group. Britton's lemma tells us that  
either the letter $b_k$ does not occur in $W$ or else $W$ 
contains an unreduced word with respect to the stable letter $b_{k}$,  
namely: 
\begin{enumerate}
\item either a sub-word of the form $b_kwb_k^{-1}$, with 
$w$ a word representing an element of $\widetilde{\varphi}_{k}(G)$; 
\item or a sub-word of the form $b_k^{-1}wb_k$, with 
$w$ a word representing an element of $i(G)$.  
\end{enumerate}
Thus $w$ is an $A$-word i.e. a word using only the letters $A$   
(constrained to belong to the image of $\phi_k$ in case (1)). 
In the first situation we choose an $A$-word $z$ so that 
$\phi_k(z)$ represents the same element in $G$ as $w$. 
If $\alpha$ is an edge loop representing a generator of $\pi_1(C)$ 
and $b_kwb_k^{-1}$ is an unreduced sub-word of $\alpha$ 
(with respect to the  HNN structure with stable letter $b_k$) 
then write $\alpha$ as $\tau w\delta$. Add the type 1 generator 
$\tau \phi_k(z) w^{-1}\tau^{-1}$ 
to the list of generators of $\pi_1(C)$ 
and replace the generator $\alpha$, with the product 
of the two generators  $\tau \phi_k(z) w^{-1}\tau^{-1}$ and 
$\alpha$ (which is $\tau b_k\phi_k(z)b_k^{-1}$). This change 
replaces the sub-word $b_kwb_k^{-1}$ of $\alpha$ by 
$b_k\phi_k(z)b_k^{-1}$. 
Now, in either case the occurrence of the unreduced sub-word  
can be eliminated  by adjoining conjugacy 2-cells along the paths labeled 
$z$ and respectively $w$. 
The new loop has fewer $B$-edges than the former one and we keep 
eliminating unreduced sub-words  until all occurrences of the letter $b_k$ 
are removed. The same method permits to get rid of all stable letters and 
hence to transform the loop into a composition of $A$-loops and hence 
a loop of type 1. 
\end{proof}

\begin{proposition}\label{top}
Suppose that $C_1$ satisfies the requirements of Proposition \ref{2gen}. 
Then there exists a finite complex $K$ obtained from $C_1$ 
by adjoining finitely many conjugation cells and finitely many 
$B$-edges $e_j$, each $e_j$ having one endpoint $u_j\not\in C_1$ such that 
$\pi_1(K)$ is generated by the classes of loops 
$\{v_1,\dots,v_m\}$ so that  each $v_s$ is  
freely homotopic in $K$ to an $A$-loop $\nu_s$
based at some $u_j$ and lying entirely in the layer $K\cap u_j X_0(G)$.
Moreover, each $\nu_j$ is disjoint from $C_1$. 
\end{proposition}
\begin{proof}
Given a compact $K$ we define the {\em layer}  
$K_x= K\cap x X(G)$.
We want to adjoin conjugacy cells in order to 
homotop all (type 1) generators of $\pi_1(C_1)$ into a disjoint 
union of layers. 

\vspace{0.1cm}
\noindent 
Define an {\em order}  on $L$ by setting 
$y <x$, if $y^{-1}x \in L$ is positive. 
This order extends to the set of layers, by saying that   
the (non-empty) layer $K_y$ is {\em below} the layer $K_x$ if 
$y < x$. We extend this terminology to (oriented) 
$B$-edges, by declaring them positive 
if their label is positive. 

\vspace{0.1cm}
\noindent 
The proof of this proposition 
follows along the lines of (\cite{BM2}, section 4). 
We define first an oriented graph $\Gamma=\Gamma_{C_1}$ 
whose vertex set $\Gamma^{0}$ is the set 
of  non-empty layers of $C_1$ i.e. slices intersecting non-trivially  $C_1$. 
The vertices $w_1,w_2$ of $\Gamma$ are joined by an oriented edge 
of $\Gamma$ if  there exists some positive $B$-edge joining 
the slice $w_1$ to the slice $w_2$. Notice that we don't ask that 
the respective layers be connected by a positive edge.  
We will consider sub-complexes $C_2$ obtained from $C_1$ 
by adding conjugacy 2-cells. Given such a complex $C_2$ and a  
subset  ${\mathcal A}$ of $\Gamma^0$ 
we say that ${\mathcal A}$ {\em carries} the loops of $C_2$ if $\pi_1(C_2)$ 
has a set of generators $\{v_1,v_2,\ldots, v_m\}$ where all 
$v_i$ are $A$-loops freely homotopic in $C_2$ to $A$-loops 
that are in the union of slices in ${\mathcal A}$.  

\vspace{0.1cm}
\noindent 
A vertex $v$ of $\Gamma$ is 
{\em extremal}  if there is no outgoing edge of $\Gamma$ issued from it. 
The key step is the following: 

\begin{lemma}
The set of extremal vertices of $\Gamma$ carries 
the loops for a suitable chosen $C_2$ which is obtained from $C_1$ by 
adjoining conjugacy cells. 
\end{lemma}
\begin{proof}
Let us start with $C_2=C_1$. Then the loops of $C_2$ are carried by 
the set of all vertices of $\Gamma$. Let ${\mathcal T}$ be a 
maximal sub-tree of $\Gamma$.  
A vertex of ${\mathcal T}$ is ${\mathcal T}$-{\em extremal} if all its 
adjacent edges  in ${\mathcal T}$ are incoming. 
Let $x$ be a non ${\mathcal T}$-extremal vertex of ${\mathcal T}$ and $xy$ 
be an oriented (outgoing) edge of ${\mathcal T}$ labeled by $b\in B$. 
Let $u$ be an $A$-loop contained in the slice  $x$. 
We add conjugacy 2-cells to $C_2$ along the $A$-loop  $u$ 
in the direction given by $b$.  In the new complex, still called 
$C_2$,  we can  freely homotop $u$ to an $A$-loop 
in the slice $y$. Notice that the slice $x$ might intersect $C_2$ 
in a non-connected sub-complex. 
By adjoining conjugacy 2-cells one might create 
additional type 1 generators, but according to Lemma \ref{add} 
the new loops are $A$-loops in the the slice $y$. 
Proceed in the same way for each non ${\mathcal T}$-extremal vertex 
of ${\mathcal T}$.  We obtain that $A$-loops in $C_2$ are carried 
by the subset of ${\mathcal T}$-extremal vertices of 
any maximal sub-tree ${\mathcal T}$.

\vspace{0.1cm}
\noindent Assume first that there is only one stable letter $b$. 
Then any maximal sub-tree ${\mathcal T}$ (and actually the graph $\Gamma$) 
is  an oriented chain since otherwise we would have a vertex of it 
with two incoming (or outgoing) $B$-edges, which is impossible. 
Moreover the  terminal vertex of the chain is both ${\mathcal T}$-extremal 
and the unique extremal vertex. The claim follows. 

\vspace{0.1cm}
\noindent
In general we can have several ${\mathcal T}$-extremal vertices, 
which might  not be extremal. 
We will show that we can get rid of those vertices which are 
${\mathcal T}$-extremal but not extremal, 
by changing the sub-tree ${\mathcal T}$ and adjoining conjugacy 2-cells. 
Let $x$ be such a vertex of ${\mathcal T}$. By hypothesis there exists 
some positive $B$-edge $e$ joining $x$ to a vertex $y$ of $\Gamma$. 
Since ${\mathcal T}$ was a maximal sub-tree the graph 
${\mathcal T}\cup e$ admits a (non-oriented) minimal length 
circuit which passes through $e$.  
Let $f\subset {\mathcal T}$ 
be the (incoming) edge of this circuit adjacent to the 
vertex $x$ and distinct from $e$. Consider then the new maximal 
sub-tree ${\mathcal T}'={\mathcal T}\cup \{e\} \setminus \{f\}$. 
The vertex $x$ is not anymore ${\mathcal T}'$-extremal and 
the set ${\mathcal A}$ is replaced by 
${\mathcal A}'={\mathcal A}-\{x\}\cup\{y\}$.  
If $y$ is ${\mathcal T}'$-extremal 
but not extremal then we continue this process and get the sequences  
$x_n$ of vertices of $\Gamma$, and ${\mathcal T}^{n}$ of maximal sub-trees.  
At some point we will find a vertex which is: 

\begin{enumerate}
\item  either non ${\mathcal T}^{n}$-extremal, and then we can reduce the 
number of ${\mathcal T}^{n}$-extremal vertices i.e. the size of the 
set carrying the loops of $C_2$; 
\item or an extremal  vertex, and we are done; 
\item or else we turn back to a vertex which has been considered 
before during this process. In this case this means the sequence of 
vertices that we meet contains an oriented circuit made 
of $B$-edges, which is impossible since $(L, <)$ is ordered. 
\end{enumerate}

\vspace{0.1cm}
\noindent
This proves that extremal vertices of $\Gamma$ carry the loops of $C_2$, as 
claimed. 
\end{proof}
 
\vspace{0.1cm}
\noindent
A layer is said to be {\em extremal} if it lies in a slice corresponding 
to an extremal vertex of $\Gamma$. 
Since $C_2$ is connected it can be arranged so that 
all extremal layers are connected. 
In fact, two points of an extremal layer can be joined by a path 
contained in $C_2$. By adding conjugacy cells we transform this path 
into an $A$-path followed by a $B$-path (or vice-versa). The $B$-path should 
be a loop and hence homotopically trivial in its own image 
since $L$ is a free group. Thus the two points are joined by a path contained 
in the layer.  Although adding 
conjugacy cells can create new loops, these are 
contained in the same connected component of the 
extremal layer. 

\begin{remark}
If there is only one stable letter then the extremal layer 
provided  by the lemma above is unique and so Proposition \ref{top} follows. 
\end{remark}

\vspace{0.1cm}
\noindent
We can suppose, without loss of generality, that 
$C_2$ contains the identity element of the group.
Recall that $p:\overline{G}\to L$ denotes 
the natural group epimorphism. 
We will need further the following technical lemma:

\begin{lemma}\label{word}
If ${C_2}_x$ is an extremal layer then we can write 
$p(x)=b y$, where $b\in B$ is a positive generator of $L$ and $y$ is an 
element of $L$ represented by a  reduced word in $B$ 
not starting with $b^{-1}$. Moreover, the layer associated to $y$ is  
not empty.  
\end{lemma}  
\begin{proof}
We have a cellular map 
$C(\overline{G},P_{\overline{G}})\to C(L,P_L)$ 
between the respective Cayley complexes induced by $p$. 
Here $P_L$ is the presentation of $L$ 
induced from $\overline{G}$. Further, we have a cellular map  
$C(L,P_L)\to C(L)$,  where $C(L)$ is the tree associated to the 
presentation $L=\langle \Phi \rangle$ with empty set of relations. 
We denote by $p:C(\overline{G},P_{\overline{G}})\to C(L)$ the composition 
of the two cellular projections.   
Observe then that the layer ${C_2}_x$ is below the layer ${C_2}_y$ 
if and only if there is a positive path from $p(x)$ to $p(y)$ in $C(L)$. 

\vspace{0.1cm}
\noindent
Since $C_2$ is connected its image $p(C_2)$ is also connected in the 
tree $C(L)$. Thus, for any $x\in C_2$ the geodesic in $C(L)$
joining the origin $1$ to $p(x)$ is contained in $p(C_2)$. We can suppose that 
the distance from $p(x)$ to the origin is at least $1$. 
Let $y$ be the vertex of this geodesic at distance $1$ from $p(x)$.  
If ${C_2}_x$ is an extremal layer of $C_2$, then 
the oriented edge $yp(x)$ of $C(L)$ is positive. 
In fact, $y\in p(C_2)$ and thus the layer ${C_2}_y$ is non-empty. 
If the edge $yp(x)$ were negative then the layer ${C_2}_x$ would be 
below the layer ${C_2}_y$, contradicting the extremality of ${C_2}_x$. 
This proves that we can write 
$p(x)=b y$, where $b\in B$ is a positive generator and $y$ is a 
reduced word not starting with $b^{-1}$.  
 \end{proof}

\vspace{0.1cm}
\noindent
For each extremal layer of $C_2$ choose a vertex $w_j$ of it and 
a  positive $B$-edge $e_j=w_ju_j$ issued from $w_j$. We adjoin the edges 
$e_j$ to $C_2$ and call $K$ the new complex.

\begin{lemma}
The  layers $K_{u_j}$ are  pairwise disjoint, disjoint from $C_2$ and  
carry the loops of $\pi_1(K)$.  
\end{lemma}
\begin{proof}
First, the slice $u_jX(G)$ through $u_j$ is disjoint from $C_2$. 
If this were not true, then ${C_2}_{w_j}$ and $C_2\cap u_jX(G)$ would represent 
vertices of $\Gamma$ connected by a positive (outgoing) edge, 
contradicting the extremality of ${C_2}_{w_j}$. 

\vspace{0.1cm}
\noindent
Next, the slices $u_{j_1}X(G)$ and $u_{j_2}X(G)$ are disjoint, for 
distinct $j_1, j_2$. 
Otherwise the two slices must coincide so that 
$p(u_{j_1})=p(u_{j_2})$. Observe that $p(u_{j})=b_{m(j)}p(w_j)$, where 
$b_{m(j)}$ is the positive generator from $B$ associated to the 
$B$-edge $e_j$. 
Lemma \ref{word} shows that 
$p(w_j)=b_{n(j)} y_j$, where $b_{n(j)}\in B$ is a positive generator 
and $y_j$ is a reduced word in the $B$ letters 
not starting with $b_{n(j)}^{-1}$.
Then we have the identity 
$b_{m(j_1)}b_{n(j_1)}y_{j_1}=b_{m(j_2)}b_{n(j_2)}y_{j_2}$ in the free group $L$. 
This implies that $j_1=j_2$. 

\vspace{0.1cm}
\noindent
Eventually, $u_j$ are extremal vertices for the graph 
$\Gamma_{K}$ associated to $K$. It suffices to show that there is no 
positive $B$-edge connecting the slices $u_jX(G)$ and $u_kX(G)$. 
If such an edge, labeled $b_i$, exists then $p(u_j)=b_i^{-1}p(u_k)$. 
Recall that $p(u_j)=b_{m(j)}p(w_j)=b_{m(j)}b_{n(j)}y_j$, where 
$y_j$ is a reduced word in the $B$ letters 
not starting with $b_{n(j)}^{-1}$ and 
 $p(u_k)=b_{m(k)}p(w_j)=b_{m(k)}b_{n(k)}y_k$, where 
$y_k$ is a reduced word in the $B$ letters 
not starting with $b_{n(k)}^{-1}$. 
Then we have the following equality  
$b_{m(j)}b_{n(j)}y_j=b_i^{-1}b_{m(k)}b_{n(k)}y_k$  holding in $L$. 
This forces $i=m(k)$ and hence $b_{m(j)}p(w_j)=p(w_k)$ holds in $L$, which 
implies that there is a positive $B$-edge labeled $b_{m(j)}$ 
joining the slices through $w_j$ and through $w_k$. In particular,  
${C_2}_{x_j}$ is not an extremal layer of $K$. This contradiction shows that 
$K_{u_j}$ are extremal layers of $K$. 

\vspace{0.1cm}
\noindent
Further $\pi_1(K)$ is isomorphic to 
$\pi_1(C_2)$ (since we simply added a number of disjoint edges) 
and the loops lying in the extremal slices of $C_2$ 
can be homotopically pushed into the slices through the $u_j$.  
Thus the set of layers $K_{u_j}$ is the set of all extremal layers 
of $K$. This proves the lemma. 
\end{proof}

\vspace{0.1cm}
\noindent
Moreover, as $K$ is connected  
it can be  arranged so that $K\cap u_j X(G)=K\cap u_j X_0(G)$ 
is connected, for every $j$. Then the previous lemmas prove the proposition.  
\end{proof}

\vspace{0.1cm}
\noindent 
{\em End of the proof of Theorem \ref{main4}.}
We have to prove that an ascending HNN-extension $\overline{G}$ of a 
finitely presented group $G$ is qsf. This is a consequence of Proposition 
\ref{top}. 
If $G$ is finitely presented then each connected component $X_0(G)$ of  
$X(G)$ is simply connected as being the Cayley complex associated 
to $G$.  Therefore  the loops $\nu_j$ are null-homotopic 
in $uX(G)$ and thus in $C(\overline{G}, P_{\overline{G}})-C$. 
This implies immediately 
that the complex $C(\overline{G}, P_{\overline{G}})$ is qsf.

\vspace{0.1cm}
\noindent 
However, when $P_G$ is infinite  
this argument does not work anymore and we need  additional ingredients.

\subsection{Constructing homotopies 
using extra 2-cells from $R_{\infty}-R$.}
Consider  now a loop $l$ 
in $K_{u}\subset uX(G)$, for an extremal layer $K_u$. 
Now $X(G)$ is the disjoint union of copies of $X_0(G)$ and each 
$X_0(G)$ embeds into the simply connected Cayley complex 
$C(G,P_G)$ of $G$.  
Therefore $uX(G)$ can be embedded in the  disjoint union of 
copies of  the simply connected 
complex $C(G, P_G)$. The later complex can be viewed as having the same 
0- and 1-skeleton as $uX(G)$, the 2-cells from $uX(G)$  and 
the additional 2-cells coming from the relations in $R_{\infty}-R$.   
Moreover the loop $l$ should be contained into one connected component
of the disjoint union. 
Thus there exists 
a null-homotopy of $l$ inside the respective $C(G, P_G)$. It is then standard 
that this implies the existence of a simplicial map 
$f:D^2\to C(G,P_G)$ from the 2-disk $D^2$ suitably triangulated, 
whose restriction to the boundary is the loop $l$. 
The image $f(D^2)$ intersects only finitely many cells of 
$C(G, P_G)$ by compactness, thus there are only finitely many 
open 2-cells $e$ of $C(G, P_G)-uX(G)$  for which the inverse image 
$f^{-1}(e)$ is non-empty. Consider the set $\{e_1(l),\ldots, e_m(l)\}$ of  
the 2-cells with this property. Each such 2-cell corresponds 
to  a relation from $R_{\infty}-R$. 
Since $f$ was supposed to be simplicial, $f^{-1}(e_{i,u}(l))$ is a finite 
union of 2-cells $e_{i_j,u}$ of the triangulated $D^2$. 
Moreover, the boundary paths 
$\partial e_{i,u}(l)$ are contained in $uX(G)$. 

\vspace{0.1cm}
\noindent 
We say that a set of  $m$ loops $\{l_j\}$ is {\em null-bordant} in $X$ if 
there exists a continuous map  $\sigma$, called a {\em null-bordism}, 
from the $m$-holed 2-sphere 
to the space  $X$ such that its restrictions to the boundary circles is $\{l_j\}$.  
In particular, the union of loops $l \cup_{i} \partial e_i(l)$ 
obtained above is   null-bordant in $uX(G)$.  Thus there exists 
a map $\sigma(l)$ from the $m$-holed sphere to $uX(G)$ whose restriction 
to the boundary is $l \cup_{i} \partial e_i(l)$, and 
we write  $l\cup_{i} \partial e_i(l)= \partial \sigma(l)$. 
We will make use of this argument further on.

\vspace{0.1cm}
\noindent Recall now that $C(G,P_G)$ was supposed to be 1-tame. 
Thus $K_u\subset \widetilde{E}_u$, where the compact $\widetilde{E}_u$ has 
the property that any loop $l\subset \widetilde{E}_u$ is homotopic within 
$\widetilde{E}_u$ to a loop $l'$ lying in $ \widetilde{E}_u- K_u$  and 
which is further 
null-homotopic in  $C(G,P_G)-K_u$. Let  $E_u=\widetilde{E}_u\cap uX(G)$, 
so that  
$E_u$ can  be written as $E_u=\widetilde{E}_u-\cup_{j=1}^ke_{j,u}$, where 
$\{e_{j,u}\}_{j=1,\ldots,k}$ is a suitable finite set of 2-cells of 
(the disjoint union of copies of) $C(G,P_G)$ 
which are not 2-cells of $uX(G)$. 
Notice that $\partial e_{j,u}\subset E_u$, for any $j$. A homotopy between the 
loop $l$ and the loop $l'$ within $\widetilde{E}_u$ induces by the 
argument above a null-bordism $H(l,l')$ between $l'$ and 
$l\cup \partial e_{j,u}(l)$ within $E_u$, where the  set $\{e_{j,u}(l)\}$ is a 
suitable subset of $\{e_{j,u}\}_{j=1,\ldots,k}$. 
Furthermore a null-homotopy of $l'$ in $C(G,P_G)-K_u$
induces a null-bordism $N(l')$ of $l'\cup \partial \delta_{j,u}(l)$ within 
$uX(G)-K_u$, where $\delta_{j,u}(l)$ are  (finitely many) 
2-cells from $C(G,P_G)-K$ 
(which are not in $uX(G)$). 

\vspace{0.1cm}
\noindent 
We consider now a finite set $J_u=\{l_{j,u}\}$ of loops which are normal 
generators of 
$\pi_1(E_u)$ and let $\{\delta_{1,u},...,\delta_{N,u}\}$ be the set of all 
2-cells $\delta_{i,u}(l)$, obtained by considering all $l\in J_u$.  

\vspace{0.1cm}
\noindent 
The key point is that $\partial e_{j,u}$ are either contained in 
$E_u$ or else disjoint from $K_u$ (and hence from $K$) while 
$\partial \delta_{j,u}$ are disjoint from $K$. 
Notice that it is the 1-tameness of  $C(G,P_G)$ 
which permitted us to discard the  
2-cells of $C(G,P_G)$ whose boundaries are not contained in $E_u$ 
but intersect $K_u$.

\subsection{Standard null-homotopies}\label{stand}
The boundary paths 
$\partial e_{i,u}, \partial \delta_{i,u} \subset uX(G)\subset C(\overline{G}, P_{\overline G})$ 
are null-homotopic within  $C(\overline{G}, P_{\overline G})$ and thus bound 
2-disks $D(e_{i,u}), D(\delta_{i,u})\subset C(\overline{G}, P_{\overline G})$.  
However  there exist some special null-homotopies for them, which are 
canonical, up to the choice of a base-point. At this place we will make use 
of the fact that the presentation $P_G$ is 
an endomorphic presentation.

\vspace{0.1cm}
\noindent Consider $\{\lambda_i\}$ be the set of loops of the 
form $\partial e_{i,u}$ or $\partial \delta_{j,u}$, for unifying the notations 
in the construction below.  
The loops $\lambda_j$  represent words which are relations 
from $P_G$ and thus can be written in the form 
\[\lambda_j =\varphi_{j_1} \varphi_{j_2}\cdots \varphi_{j_{k_j}}(r_{\alpha_j}) \] 
where $r_{\alpha_j}\in R$ and the $j_i$'s depend on $j$. 
We have implicitly chose the convenient orientation of the loops $\lambda_j$ 
in order to be recovered from $r_{\alpha_j}$ and not from 
$r_{\alpha_j}^{-1}$.  
It is important to notice that all $\varphi_j$ appear only with positive exponents in the expression above. 
Recall that $R$ is the set of relations 
that survive in $P_{\overline{G}}$.  We can identify a loop with the word that 
represents that loop in the Cayley complex. Thus, by abuse of notation,  
we can speak of $\varphi_k(l)$ where $l$ is a loop. 
Observe that the loop 
$\varphi_m(l)$ is freely homotopic to the loop $l$, since it is associated  
to a specific conjugate  in terms of words. This homotopy is  
the cylinder $C_m(l)$  which is the union of all conjugacy cells based 
on elements of $l$ and using the vertical element $b_m$. 
The loop corresponding  to 
$\varphi_{j_1} \varphi_{j_2}\cdots \varphi_{j_{k_j}}(r_{\alpha_j})$ 
is one boundary of the cylinder $C_{j_1}(\lambda_j^{(1)})$. The  other 
boundary of this cylinder is the loop 
$\lambda^{(1)}_j=\varphi_{j_2}\cdots \varphi_{j_{k_j}}(r_{\alpha_j})$. 
The second loop has, in some sense, smaller complexity  
than the former one and we  can continue to simplify it.     
The cylinder $C_{j_s}(\lambda_j^{s})$ interpolates between  
$\lambda_j^{(s)}=\varphi_{j_{s+1}}\cdots \varphi_{j_{k_j}}(r_{\alpha_j})$ 
and $\lambda_j^{(s-1)}$. Set 
$C(\lambda_j)=\cup_{1\leq s\leq k_j} C_{j_s}(\lambda_j^{(s)})$. 
Eventually, recall that 
$r_{\alpha_j}\in R$ and thus the corresponding loop bounds a  2-cell 
$\varepsilon_{\alpha_j}$ of $X(G)$. Thus $D(\lambda_j)=C(\lambda_j) 
\cup \varepsilon_{\alpha_j}$ 
is a specific 2-disk giving an explicit null-homotopy of $\lambda_j$ 
within  
$C(\overline{G},P_{\overline{G}})$.

\subsection{Saturation of layers}
\vspace{0.1cm}
\noindent 
Given a compact $K$ we considered the layers $K_x=K\cap x X(G)$. 
Observe, following \cite{BM1}, that we can suppose that 
all intersections $K\cap x X_0(G)$ are connected for all $x$ where non-empty,
and $K_x\cap xX_0(G)=K\cap xX_0(G)$ if $K_x$ is an extremal layer.

\vspace{0.1cm}
\noindent 
The finite complex $K$ is said to be {\em saturated} if it has the 
following property. For each vertex $c$ of $C$ and 
positive $B$-path at $c$ that  ends at $c'$ in an extremal layer of $K$ 
the endpoint $c'$ is in $K$.



\begin{lemma}
We can assume that the complex $K$ obtained in Proposition \ref{top}
has saturated layers. 
\end{lemma}
\begin{proof}
It suffices to add finitely many conjugacy cells in order to achieve 
the saturation. Moreover, when adjoining conjugacy cells we do not create 
extra loops of type 2 and hence the requirements in 
Proposition \ref{top} are still satisfied. 
\end{proof}

\vspace{0.1cm}
\noindent
Recall now that $\partial \delta_{j,u}$ are disjoint from $K_u$, for 
any extremal layer $K_u$. 
We have then: 

\begin{lemma}\label{sat}
If $K$ is saturated then $D(\partial \delta_{j,u}) \cap C =\emptyset$.
\end{lemma} 
\begin{proof}
If we had a point $c$ belonging to $D(\partial \delta_{j,u}) \cap C$ 
then there would  exist a path from $c$ to a point $c'$ 
in $\partial \delta_{j,u}$, which contains 
only vertical segments from  the cylinder 
$C(\partial \delta_{j,u})$. This is then a positive $B$-path  
and thus its endpoint $c'$ belongs to $K_u$, by the saturation hypothesis, 
but  this contradicts the fact that   $\partial \delta_{j,u} \subset E_u-K_u$.  
\end{proof}

\begin{remark}
The analogous statement fails in the case when we take for $R_{\infty}$ 
the set of {\em reduced} words in the free group $F(A)$ coming from 
iterating the $\phi_i$ on the set $R$, in general.   
\end{remark}

\vspace{0.1cm}
\noindent Although $\partial e_{j,u}$ might intersect $K$ they are contained in 
$E_u$. Consider then 
$W=\{(\alpha,u); D(\partial e_{\alpha,u}) \cap C \neq \emptyset\}$. 
We construct therefore the following set: 
\[ Z =K\cup_{u} E_u \cup_{(\alpha,u) \in W} D(\partial e_{\alpha,u}) \]

\begin{lemma}\label{pi}
If $K$ is saturated then the inclusion $K\cup_{u}E_u \hookrightarrow Z$ 
induces a surjection $\pi_1(K\cup_{u}E_u) \to \pi_1(Z)$. 
\end{lemma}
\begin{proof}
The only new loops appearing when we added the cylinders $C(\partial 
e_{\alpha,u})$ come either 
from their intersections with $K$ or else from their pairwise intersections. 

\vspace{0.1cm}
\noindent 
In the first case consider $q\in D(\partial e_{\alpha,u}) \cap K\neq \emptyset$. 
The new loop $\lambda$ created this way is the composition of 
an $A$-path from a vertex $*$ of $K_u$ to a vertex of $\partial e_{\alpha,u}$ 
followed by a $B$-path in $C(\partial e_{\alpha,u})$ and then by 
a path in $K$ to the point $*$. 
Now $(\alpha,u)\in W$, so  that 
$C(\partial e_{\alpha,u}) \cap C \neq \emptyset$. 
Any vertex  of $D(\partial e_{\alpha,u}) \cap K$ belongs therefore 
to a positive $B$-path starting 
at a point of $C$ and ending at the extremal layer $K_u$.
Thus we can homotopically push such a loop $\lambda$ using the 
conjugacy cells    --  that are  
contained both in the cylinders $C(\partial e_{\alpha,u})$ and in $K$, 
since $K$ is saturated -- until they reach the layer $Z\cap u X(G)=E_u$. Thus 
the subgroup generated by images of $A$-loops in $K$ and loops in $E$ 
contains the loops of the form $\lambda$ from $\pi_1(Z)$.

\vspace{0.1cm}
\noindent 
In the second case assume that $C(\partial e_{\alpha,u})\cap 
C(\partial e_{\beta,v})\neq \emptyset$. If $K_u=K_v$ 
the proof from above applies without essential modifications. 
This proves the lemma for the case when we have only one stable letter.
 
\vspace{0.1cm}
\noindent   
Assume now that $K_u\neq K_v$. 
Let $q$ be an intersection point 
of these cylinders. A loop $\lambda$ created by this double point 
is then the composition of 
an $A$-path joining a point $*$ of $K_u$ to some point of $E_u-K_u$, 
followed by a $B$-path in $C(\partial e_{\alpha,u})$ 
reaching $q$ then a $B$-path in $C(\partial e_{\beta,v})$ to a point 
in $E_v-K_v$, followed by an $A$-path to a point of $K_v$ and eventually 
by a path in $K$ joining it to $*$. 
The only problem, with respect to the previous case, is that we cannot 
push directly the loop $\lambda$ along conjugacy cells since 
we have two extremal layers. The idea is to decompose it as the 
composition of two loops, each one of them which can be 
homotopically pushed into one extremal layer. 
 
\vspace{0.1cm}
\noindent 
The subset $p(K)\subset C(L)$ is connected and thus the geodesic 
$\gamma$ in $C(L)$ joining $p(u)$ and $p(v)$ is contained in $p(K)$. 
The cylinder $C(\partial e_{\alpha,u})$ (respectively 
$C(\partial e_{\beta,v})$) is made of conjugacy cells starting 
from some relation in $R$ 
in the direction given by a positive $B$-path $\gamma_1$ 
(respectively $\gamma_2$), as explained in section \ref{stand}. 
Then $\gamma_j$ are positive paths and hence geodesics in the tree $C(L)$
having  a common vertex, namely $p(q)$. Therefore  
$\gamma$, $\gamma_1$ and $\gamma_2$ have a common vertex, say $y$. 
Further the positive $B$-sub-paths $\gamma_1[p(q)y]$ and $\gamma_2[p(q)y]$ 
coincide. We can push 
$q$ along this positive $B$-sub-path and find that  
$C(\partial e_{\alpha,u})\cap C(\partial e_{\beta,v})$ 
contains also a vertex $t$ in the slice associated to $y$, namely 
with $p(t)=y$. It suffices then to consider the case where $p(q)=y$.  

\vspace{0.1cm}
\noindent
Furthermore, we know that $y\in p(K)$, which implies that the layer 
$K_t$ is not empty. We claim that: 
\begin{lemma}\label{exist}
There exists a vertex of $z\in K_t$ 
which is in the same connected component of $tX(G)$ as $t$. 
\end{lemma}
\begin{proof}
Let $t_u\in E_u$ (respectively $t_v\in E_v$) be the endpoint of the 
$B$-path given by the word $\gamma_1[yp(u)]$  (respectively 
$\gamma_2[yp(v)]$) starting at $t$.  Recall now that 
$E_u$ and $E_v$ are each connected and thus we can find vertices 
$w_u\in K_u$, $w_v\in K_v$  which are joined to $t_u$ and $t_v$ respectively 
by $A$-paths  corresponding to words $a_u$ and $a_v$ in the $A$-letters. 

\vspace{0.1cm}
\noindent
Observe that  the $B$-sub-paths 
$\gamma[p(u)y]$ and $\gamma_1[p(u)y]$ joining $p(u)$ and $y$ 
in $C(L)$ coincide, as well as 
$\gamma[p(v)y]$ and $\gamma_2[p(v)y]$. 

\vspace{0.1cm}
\noindent
Consider then a path joining $w_u$ to $w_v$ in the connected 
sub-complex $p^{-1}(\gamma)\cap K$. This path is given by a word of the 
following form: 

\[ U=a_{2k+1}b_{i_{2k}} a_{2k}b_{i_{2k-1}}\cdots 
 a_{k+2}b_{i_{k+1}}a_{k+1}b_{i_{k}}^{-1}a_{k}b_{i_{k-1}}^{-1}\cdots 
a_{2}b_{i_{1}}^{-1}a_{1}\]
where $a_j$ are words in the $A$-letters and $b_i\in B$ are 
the positive generators. Furthermore the $B$-path 
$\gamma[yp(v)]$ is given by the word 
$B_+= b_{i_{2k}}b_{i_{2k-1}}\cdots b_{i_{k+1}}$, while 
the $B$-path $\gamma[p(u)y]$ is given by the word 
$B_-=b_{i_{k}}^{-1}b_{i_{k-1}}^{-1}\cdots b_{i_{1}}^{-1}$. 
Notice that $b_{i_{k+1}}\neq b_{i_k}$ since $\gamma$ is a geodesic. 
We have then a loop $\lambda_0$ in the Cayley graph of $\overline{G}$ 
realizing the word $B_+^{-1}a_v^{-1}Ua_u^{-1}B_-^{-1}$. 
This word must therefore represent the identity in $\overline{G}$. 
We use induction on $k$ and Britton's lemma 
to obtain that the only way that this word can be simplified 
to the empty word is by means of reductions 
of the type $b a b^{-1}=c$, where $\phi_j(c)=a$, for $b\in B$ 
and $a,c\in A$. 
This means that there exist families of 
conjugacy cells in $C(\overline{G},P_{\overline{G}})$, where the first  
family touches the extremal slice along the path 
$a_1a_u^{-1}$ (respectively $a_{2k+1}a_v^{-1}$) in  
the direction $b_{i_1}^{-1}$ (respectively $b_{i_{2k}}^{-1}$) 
and the next ones  use inductively 
the directions given by $b_{i_2}^{-1},\cdots,b_{i_k}^{-1}$ 
(respectively $b_{i_{2k-1}}^{-1},\cdots,b_{i_{k+1}}^{-1}$). 
Each family connects one slice to the slice below it. 
We can therefore push homotopically in $C(\overline{G},P_{\overline{G}})$ 
the loop $\lambda_0$ along these conjugacy cells to the lowest  slice $tX(G)$.  
But, at each step, the pushed loop has at least 
one vertex from $K$. Thus there exists 
a vertex $z\in K_v$ which is connected by an $A$-path 
to the vertex $t$, as claimed. 
\end{proof}

\vspace{0.1cm}
\noindent
We turn back now to the loop $\lambda$. Since $p(\lambda)\subset C(L)$ 
contains both $p(u)$ and $p(v)$, 
it should contain the geodesic $\gamma$ and then, by the 
previous arguments, there exists a point $z'$ of $\lambda$ in the 
layer $K_v$. Using Lemma \ref{exist} for the points $z$ and $z'$ 
instead of $z$ and $t$ it follows that $z'$ and $z$ 
(and hence $t$) are in the same connected component $tX_0(G)$ of $tX(G)$  
so that $z'$ can be joined by an $A$-path $\zeta$ to $t$. 
Therefore we can split $\lambda$ as the composition of two loops 
$\lambda_1\lambda_2$, by inserting $\zeta$ between $z'$ and $t$.
But now each one of the two loops $\lambda_i$ can be homotopically pushed 
within $K$ in the directions given by the $B$-sub-paths of $\gamma[yp(u)]$ and  
respectively $\gamma[yp(v)]$ to one of the extremal slices $E_u$ or $E_v$.  

\vspace{0.1cm}
\noindent 
This proves that $A$-paths in $K$ and loops in $\cup_uE_u$ generate 
all of $\pi_1(Z)$. This settles Lemma \ref{pi}. 
\end{proof}

\subsection{End of the proof of Theorem \ref{ext}}  Take a loop 
$l$  in $\pi_1(Z)$. It can be supposed that $l$ is either from the set 
$\sqcup_uJ_u$ 
that normally generates $\pi_1(\sqcup_uE_u)$ (recall that $E_u$ are disjoint) 
or else from the set of $A$-loops $\nu_j$ furnished by Proposition \ref{top}. 
We observed that $l\subset E_u$ and 
$l'\cup_j \partial e_{j,u}(l)\subset E_u-K_u$ are null-bordant in  
$E_u\subset uX(G)$ using the  
null-bordism  
$H(l,l')$. Moreover $l$ and $l''=l'\cup_{j; (j,u)\not\in W} e_{j,u}(l)$ are 
null-bordant in $Z$  by means of the modified null-bordism 
$H(l,l')\cup_{(j,u)\in W} D(e_{j,u}(l))$, since 
the boundaries $\partial e_{j,u}$, with $(j,u)\in W$, are null-homotopic in 
$Z$.  Moreover, $l'\cup_{j; (j,u)\not\in W} e_{j,u}(l)\cup_k \partial 
\delta_{k,u}$ 
is  furthermore null-homotopic in 
$uX(G)-K_u\subset  C(\overline{G}, P_{\overline{G}}) -C$. 
We adjoin then the 2-disks $D(e_{j,u}(l))$ and 
$D(\partial \delta_{k,u})$ and obtain 
a null-homotopy of $l''$ within $C(\overline{G}, P_{\overline{G}}) -C$. 
This means that $C(\overline{G}, P_{\overline{G}})$ is 1-tame and thus 
qsf. This proves Theorem \ref{ext}. 

\begin{remark}
\begin{enumerate}
\item We can always add new relations to  the group presentation $P_G$ 
in order to make it 1-tame. However the new presentation 
is not necessarily a finite endomorphic presentation.  
Thus the  second assumption in Theorem \ref{ext} seems nontrivial. 
\item We don't know whether the 1-tameness of a presentation $P$ with 
infinitely many relations  which are unreduced words  
is equivalent to the 1-tameness 
of the presentation $P_r$ consisting of the reduced words 
arising in the relations of $P$. It does so, for instance, 
when the length of the cancelling sub-words (i.e. sub-words
of the form $a_1a_2\cdots a_ka_k^{-1}\cdots a_2^{-1}a_1^{-1}$)
is uniformly bounded. 
\end{enumerate}   
\end{remark}

\section{Proof of Theorem \ref{main3}}
\subsection{The Grigorchuk group. } 
 Grigorchuk constructed  in the eighties
a finitely generated infinite torsion group of intermediate growth
having solvable word problem (see \cite{Gri1}). This group is not
finitely presented but Lys\"enok  obtained  (\cite{Ly}) a nice
recursive presentation of $G$ as follows:
\[ G=\langle a,c,d \;|\; \sigma^n(a^2), \,\sigma^n((ad)^4),\,
\sigma^n((adacac)^4), \, n\geq 0 \rangle \] where
$\sigma:\{a,c,d\}^*\to \{a,c,d\}^*$  is the substitution  that
transforms words according to the rules: 
\[ \sigma(a)=aca, \sigma(c)=cd, \sigma(d)=c\]
We denote by $A^*$ the set of positive nontrivial words in the letters 
of the alphabet $A$ i.e. the free monoid generated by $A$ without the 
trivial element. 

\vspace{0.1cm}\noindent 
The finitely  presented  HNN-extension $\overline{G}$ of the Grigorchuk 
group $G$   was  constructed for
the  first  time in \cite{Gri0,Gri} as  a  group  with  5  
generators  and 12  short
relations. The group $\overline{G}$ is a finitely presented 
example of a group which is amenable but not elementary amenable.
Bartholdi  transformed  this  presentation  in  
the  form of a presentation
with 2  generators  and  5 relations, as described in 
\cite{DCG}. Later Bartholdi  presented (see \cite{Bar})  
some  general  method of  getting endomorphic  
presentations  for  branch  groups.

\vspace{0.1cm}\noindent 
The endomorphism of $G$ defining the  HNN extension is induced by the
substitution $\sigma$ and thus the new group $\overline{G}$ has the following
finite presentation:
\[ \overline G = \langle a,c,d,t \; | \;
a^2=c^2=d^2=(ad)^4=(adacac)^4=1, \; a^t=aca, c^t=dc, d^t =c
\rangle \]
where $x^y=yxy^{-1}$. 
Theorem \ref{ext} is the main ingredient needed for proving 
Theorem \ref{main3}, which we restate here for the sake of completeness: 

\begin{theorem}\label{gri}
The HNN extension of the Grigorchuk group is qsf. 
\end{theorem}

\begin{remark}
Relations in the Lys\"enok endomorphic presentation of   
Grigorchuk's group  are given by iterating the substitution $\sigma$ 
and thus  involve only words with positive exponents 
on the generators which are reduced words.  
\end{remark}

\subsection{The  Lysen\"ok presentation is 1-tame}
We want to use Theorem \ref{ext}. Using the notations from section 5 
the group $L$ is the infinite cyclic group generated by 
the endomorphism $\sigma$.  
Since the endomorphism $\sigma$ is expansive 
there are only finitely many positive paths between two elements 
of $L$. Further, the map $M\to L$ is obviously injective.  

\vspace{0.2cm}
\noindent
In the next section we will show that: 
\begin{proposition}\label{infgri}
The group $G$ with the Lys\"enok presentation $P_G$ is 1-tame. 
\end{proposition}

\vspace{0.2cm}
\noindent 
This proposition and Theorem \ref{ext} will settle then  
the proof of Theorem \ref{main3}. 

\noindent The main idea is that the group $G$ is commensurable with 
$G\times G$ (see e.g. \cite{DLH}, VIII.C. Theorem 28, p.229). 
Further, the qsf property is invariant under commensurability. 
Moreover, the proof in \cite{BM2} which shows that extensions
of infinite finitely presented groups are qsf works also 
in the infinitely presented case. Even more, \cite{BM2} shows
that extensions of infinite finitely presented groups are actually 1-tame. 
Thus $G\times G$ with any 
direct product presentation  is 1-tame. In particular, this happens 
if we consider the presentation $P_{G\times G}=P_G\times P_G$, defined 
as follows: 
\begin{itemize}
\item take two copies of the generators, $a_j,b_j,c_j,d_j$, $j\in\{ 1,2\}$, 
corresponding to $G\times \{1\}$ and $\{1\}\times G$ respectively; 
\item take two copies of the Lys\"enok relations corresponding to 
each group of generators.  
\item add the commutativity relations between generators from distinct groups, 
namely: 
\[ [a_1,a_2]=[a_1,b_2]=[a_1,c_2]=[a_1,d_2]=1\]
and the similar ones involving $b_1,c_1$ and $d_1$. 
\end{itemize}
\vspace{0.2cm}
\noindent Since $G$ is commensurable to $G\times G$ we will show that 
the presentation $P_{G\times G}$ induces a presentation $P_G^*$ of $G$. 
The induction procedure consists of transferring presentations 
towards - or from -  a finite index normal subgroup and transport it 
by some isomorphism.  
In particular, $G$ with the induced presentation $P^*_G$ is 1-tame.  

\vspace{0.2cm}
\noindent
We will show below that the 
$P^*_G$ (up to finitely many relations)  consists of $P_G$ 
and finitely many families of relations, each family being  
conjugated to the family of standard relations 
in $P_G$. 
The later relations can be simply discarded from $P_G^*$ 
without affecting the qsf property of the associated Cayley complex. 
In particular, the presentation obtained after that is in the 
same finite equivalence class as $P_G$.  
This will imply that the group  $(G,P_G)$ is qsf and thus 
its HNN extension $\overline{G}$ is also qsf, according 
to the Theorem \ref{ext}.

\begin{remark}
Other examples of groups with endomorphic presentations  
including branch groups are given in \cite{Bar}. Our present methods 
do not permit handling all of them. It is very probable that a general 
method working for this family will actually yield the fact that 
any finitely presented group admitting a normal (infinite) finitely generated 
subgroup of infinite index should be qsf. 
\end{remark}

\subsection{Preliminaries concerning $G$ following (\cite{DLH},VIII.B)}
It is customary to use the following 4 generators presentation of $G$: 
\[P_G(a,b,c,d):  G=\langle a,b,c,d \;|\; a^2=b^2=c^2=d^2=bcd=1, \,\sigma^n((ad)^4),\,
\sigma^n((adacac)^4), \, n\geq 0 \rangle \] where
$\sigma:\{a,b,c,d\}^*\to \{a,b,c,d\}^*$  is the substitution  that
transforms words according to the rules:
\[ \sigma(a)=aca,\sigma(b)=d, \sigma(c)=b, \sigma(d)=c\] 
from which we can drop either the generator $b$ or else $c$ and get 
the equivalent presentations with three generators: 
\[P_G(a,c,d):  G=\langle a,b,c,d \;|\; a^2=c^2=d^2=1, \,\sigma^n((ad)^4),\,
\sigma^n((adacac)^4), \, n\geq 0 \rangle \] where
$\sigma:\{a,c,d\}^*\to \{a,c,d\}^*$  is the substitution  that
transforms words according to the rules:
\[ \sigma(a)=aca,\sigma(c)=cd, \sigma(d)=c\] 
or else: 
\[ P_G(a,b,d):  G=\langle a,b,d \;|\; a^2=b^2=d^2=1, \,\sigma^n((ad)^4),\,
\sigma^n((adabdabd)^4), \, n\geq 0 \rangle \] where
$\sigma:\{a,b,d\}^*\to \{a,b,d\}^*$  is the substitution  that
transforms words according to the rules:
\[ \sigma(a)=abda,\sigma(b)=d, \sigma(d)=bd\] 

\vspace{0.1cm}
\noindent Define $G^0$ be the subgroup consisting of words in $a,b,c,d$ 
having an even number of occurrences of the letter $a$. This is the same 
as the subgroup denoted $St_G(1)$ in (\cite{DLH}, VIII.B.13 p.221). It is 
clear that $G^0\lhd G$ is a normal subgroup and we have an exact sequence: 
\[ 1\to G^0 \to G \to G/G^0= \Z/2=\langle a \rangle  \to 1\]
where $G/G^0$ is generated by $aG^0$. 
It follows that $G^0$ is the subgroup of $G$ generated by 
the following 6 elements: 
\[ G^0=\langle b,c,d, aba, aca, ada \rangle \subset G\]

\vspace{0.1cm}
\noindent There exists an injective homomorphism 
$\psi: G^0\to G\times G$ given by the formulas:
\[ \psi(b)=(a,c), \psi(c)=(a,d), \psi(d)=(1,b)\]
\[ \psi(aba)=(c,a), \psi(aca)=(d,a), \psi(ada)=(b,1)\]

\vspace{0.1cm}
\noindent
Let $B\lhd G$ be the normal subgroup generated by $b$. It is known 
that $B=\langle b, aba, (bada)^2, (abad)^2\rangle$. We have then an exact 
sequence: 
\[ 1 \to B \to G \to G/B=D_8=\langle a, d\rangle \to 1 \]
where $G/B$ is the dihedral group of order 8, denoted $D_8$. Moreover $D_8$ is 
generated by the images of the generators $a$ and $d$. Since the subgroup 
$D=\langle a, d\rangle\subset G$ is the dihedral group $D_8$ it actually 
follows that the extension above is split. 
Consider further 
the group $D^{diag}=\langle (a,d), (d,a)\rangle \subset G\times G$ 
which is isomorphic to the group $D_8$. 
Then we can describe the image of $\psi$ as  
$\psi(G^0)=(B\times B)\ltimes D^{diag}\subset G$. Notice that the later is 
a subgroup (although not a normal subgroup) of $G\times G$ having index 8.

\vspace{0.1cm}
\noindent It is easier to work with normal subgroups below since we 
want to track explicit presentations in the commensurability process. 
Therefore we will be interested in the subgroup 
$B\times B\subset \psi(G^0)\subset G\times G$ which is a normal subgroup. 
Denote by $A$ the inverse image $\psi^{-1}(B\times B)$ which is a normal 
subgroup of $G^0$. It follows that $G^0/A \to 
(B\times B)\ltimes D^{diag}/B\times B=D$ is an isomorphism and 
$G^0/A$ is generated by the images of $c$ and $aca$. Moreover the subgroup 
$\langle c, aca\rangle  \subset G^0$ is dihedral of order 8 and thus 
there is a split exact sequence: 
\[ 1 \to A \to G^0 \to G^0/A=D_8=\langle c, aca \rangle \to 1\]
Collecting these facts it follows that actually $A$ is the 
normal subgroup of $G$ generated by $d$ and we have a split exact 
sequence: 
\[ 1 \to A \to G \to G/A=D_{16}=\langle a,c\rangle\to 1\]
where $G/A$ is generated by the images of $a$ and $c$ and it is isomorphic 
to the group $E=\langle a,c\rangle \subset G$, which is the dihedral group of order 16.

\subsection{Inducing group presentations} The presentation $P_{G\times G}$ 
of $G\times G$ induces a presentation $P_{B\times B}$ 
of its normal finite index subgroup $B\times B$. 
The isomorphism $\psi:A \to B\times B$ 
transports $P_{B\times B}$ into the presentation $P_A$ of $A$. Eventually 
we can recover the presentation $P^*_G$ of $G$ 
from that of its normal subgroup $A$. 
In order to proceed we need to know how to induce  
presentations from and to normal finite index subgroups. 

\vspace{0.1cm}
\noindent First we have the following well-known lemma of Hall: 
\begin{lemma}  
Assume that we have an exact sequence: 
\[ 1 \to K \to G \to F\to 1 \]
and $K=\langle k_i | R_j \rangle$, 
$F=\langle m_j | S_n \rangle$  are group presentations in the generators 
$k_i$ (respectively $m_j$) and relations $R_j$ (respectively $S_n$).  
Then $G$ has a presentation of the following form: 
\[ G = \langle k_i, m_j| R_j, S_n(m_j)=A_n(k_i), m_j k_i m_j^{-1}= 
B_{ji}(k_i)\rangle \]
where $A_n, B_{ji}$ are suitable words in the generators $k_i$. 
Specifically the relations using $A_n$ express the relations between 
the lifts of the generators $m_j$ to $G$, while the last relations 
express the normality of $K$ within $G$. 
\end{lemma}

\vspace{0.1cm}
\noindent Inducing presentations to a normal subgroup seems slightly 
more complicated. For the sake of simplicity we formulate the 
answer in the case where the relations are positive 
(i.e. there are no negative exponents) and the extension is split, as it 
is needed for our purposes. Observe however that the result can be extended 
to the general situation. 

\begin{lemma}  
Assume that we have a split exact sequence: 
\[ 1 \to K \to G \to F\to 1 \]
where $G=\langle \{x_i\}_{i=1,...,N} | R_j \rangle$ and the group 
$F$ is finite.  
Let $F=\{1=f_0,f_1,f_2,\ldots, f_n\}$ be an enumeration of its elements. 
Assume further that the projection map $p:G\to F$ takes the 
form $p(x_i)=f_{p(i)}$ where $p$ is a map 
$p:\{1,2,\ldots, N\}\to \{0,1,\ldots, n\}$. 
Assume that the relations $R_j$ read as 
\[ R_j= x_{i_1}\cdot x_{i_2}\cdots x_{i_{k+1}}\]
We choose lifts $\widehat{f_j}\in G$ for the elements $f_j$, using the 
splitting homomorphism.  
Set then $y_j=x_j \widehat{f_p(j)}^{-1}$ and denote by 
$^{f_k}y_j=\widehat{f_k}y_j\widehat{f_k}^{-1}$ the conjugation. 
We consider below  \,$^{f_k}y_j$ as being distinct symbols, called 
$y$-letters, for all $k\neq 0$ and $j$. 

\vspace{0.1cm}
\noindent 
Then the group $K$ admits the following presentation: 
\begin{itemize} 
\item The generating set is the set of $N(n+1)$ elements 
$y_j$, $^{f_k}y_j$, $k\in \{1,2,\ldots n\},  j\in\{1,2,\ldots, N\}$. 
\item Relations are obtained using the following procedure. 
\begin{itemize}
\item Each relation  $R_j= x_{i_1}\cdot x_{i_2}\cdots x_{i_k}$ gives rise  
to a basic relation in the $y$-letters alphabet: 
\[ R'_j= y_{i1}\cdot \,\left(^{f_{i_1}}y_{i_2}\right)\cdot \,\left(^{f_{i_1}\cdot f_{i_2}}y_{i_2}\right)\cdots \,\,\left(
^{f_{i_1}\cdot f_{i_2}\cdots f_{i_k}}y_{i_k}\right)\]
where each superscript product $f_{i_1}\cdot f_{i_2}\cdots f_{i_s}$ 
is replaced by its value, as an element  
$f_{\lambda(i_1,i_2,...,i_s)}\in F$. 
\item Next one considers all images of the basic relations $R'_j$ under 
the action of $F$ (by conjugacy).  
Specifically, for any basic relation in $y$-letters 
\[ R= ^{f_{j_1}}y_{j_1}\cdot\, ^{f_{j_2}}y_{j_2}\cdots  \,^{f_{j_p}}y_{j_p}\] 
and any element $f\in F$  
one associates the relation 
\[ ^fR=  ^{ff_{j_1}}y_{j_1}\cdot\, ^{ff_{j_2}}y_{j_2}\cdots \, ^{ff_{j_p}}y_{j_p} \]
in which each superscript is considered as an element of $F$. 
\end{itemize}
\end{itemize} 
\end{lemma}

\noindent Here we set the notation $^ay$ in order to emphasize that these 
are abstract symbols, which will be viewed as elements of $K$. 
They will be equal to the usual  conjugacies $y^a$ only when 
seen as elements of $G$.  
\begin{proof}
Any element of $G$ is a product of $y$-elements and some $\widehat{f_j}$. 
Thus an element of $K$ should involve no $\widehat{f_j}$. 

\vspace{0.1cm}
\noindent 
Remark now that expressing $R_j$ using the elements \,$^{f_k}y_j$ we obtain 
\[ R_j= \left (y_{i1}\cdot \,^{f_{i_1}}y_{i_2}\cdot \,^{f_{i_1}\cdot f_{i_2}}y_{i_2}\cdots \,\,
^{f_{i_1}\cdot f_{i_2}\cdots f_{i_k}}y_{i_k}\right) 
\widehat{f_{i_1}}\cdot \widehat{f_{i_2}}\cdots \widehat{f_{i_k}}\]
Moreover, the product of the first $k$ terms in the right hand side 
is an element of $K$. Since the extension is split we should have 
$\widehat{f_{i_1}}\cdot \widehat{f_{i_2}}\cdots \widehat{f_{i_k}}=1$ 
coming as a relation in $F$. Thus $R'_j=1$, as claimed. 
It is clear then that $^fR'_j=1$ holds true also because 
$K$ is a normal subgroup. 

\vspace{0.1cm}
\noindent 
In order to see that these relations define $K$, consider the 2-complex
$Y_G$ associated to the given presentation of $G$. Thus $Y_G$ has 
one vertex $v$.  Then  
$K$ is the fundamental group of the finite covering $\widehat{Y_G}$ 
(with deck group $F$) of $Y_G$, that is associated to the projection map 
$G\to F$. This is a non-ramified covering of degree $|F|$, the order of $F$.   
Thus each open 2-cell of $Y_G$ is covered by precisely $|F|$ 2-cells 
of  $\widehat{Y_G}$. It would suffices now to read the presentation of 
$\pi_1(\widehat{Y_G})$ on the cell structure of  $\widehat{Y_G}$. 
The only problem is that loops in $Y_G$  
lift to paths in $\widehat{Y_G}$ 
which are not closed. Now $\widehat{Y_G}$ has $|F|$ vertices 
that are permuted among themselves by $F$, let us call them 
$v^f$, for $f\in F$, such that the deck transformations 
act as $g\cdot v^f=v^{gf}$. The vertex $v^1$ will be the base point of 
$\widehat{Y_G}$. 
The loop $l_j$ based at $v$ that corresponds to the 
generator $x_j$ lifts to a path $c_j$  joining $v^1$ to $v^{p(j)}$. 
Moreover the inverse image of the loop $l_j$ under the covering is the 
union of all translated copies $fc_j$ (joining $v^f$ to $v^{fp(j)}$) 
of this path, which should be distinct as the covering is non-ramified. 
In this setting we have a natural presentation of $\pi_1(\widehat{Y_G})$ 
as a fundamental groupoid with base-points $v^f$, for all $f\in F$. 
Simply take all (oriented) edges $fc_j$ as generators and all 
2-cells as relations. The 2-cells are all disjoint and permuted among 
themselves by $F$ and in each $F$-orbit the 
2-cell based at $v^1$ corresponds to one 2-cell of $Y_G$. 
One could choose now a maximal tree (corresponding to the choice of the 
elements $\widehat{f_j}$) in the 1-skeleton of $\widehat{Y_G}$ 
and collapse it in order to find a complex which comes from a 
group presentation. 
Alternatively we can transform the groupoid presentation into a 
group presentation by choosing a fixed set of paths $l(f)$ joining 
$v^1$ to $v^f$. The choice of this system amounts to choose lifts 
$\widehat{f_j}$ in $G$. Then the paths  
$l(f)\cdot fc_j l(fp(j))^{-1}$ are now based at $v^1$ and represent 
a generator system for the loops in $\widehat{Y_G}$. This loop represents  
the generator $^fy_j$ of $K$ under the identification with  
$\pi_1(\widehat{Y_G})$. Further the 2-cell based at $v^1$ corresponds 
to the basic relation associated to a relation in $G$ and its 
images under the deck transformations are those described in the 
statement.  Thus the fundamental group $\pi_1(\widehat{Y_G})$ based 
at $v^1$ has the claimed presentation.    
\end{proof}

\subsection{Carrying on the induction for the Grigorchuk group}
We will consider first the group $G$ with its presentation $P_G(a,b,d)$ 
and the normal subgroup $B$ normally generated by $B$. According to the 
induction lemma above we have a natural system of generators given 
by $^{G/B}b=^{\langle a, d \rangle}d$ which is simply a notation for 
\[ \{^{x}b; x\in G/B =\langle a, d\rangle\}= \{b,\, ^ab,\, ^{d}b, \,
^{ad}b,\, ^{da}b,\, ^{ada}b,\, ^{dad}b, \,^{adad}b \}\]
The infinite set of words 
$w_n=\sigma^n((ad)^4),\, z_n=\sigma^n((adabdabd)^4)$ 
are relations in $G$ that induce relations $T(w_n)$ and $T(z_n)$ 
in $B$, by the procedure above. This  amounts to the following. 
Write first $w_n$ (and $z_n)$ as a word in $a,b,d$ as follows: 
\[ w_n =w_{n,0}(a,d)\, b\, w_{n,1}(a,d)\,b\,\cdots w_{n,k}(a,d)\,b\, w_{n,k+1}(a,d)\] 
where $w_{n,j}(a,d)$ are words in $a$ and $d$ and thus can be reduced 
as elements of $D=G/B$. 
Then the basic relation corresponding to $w_n$ is now 
\[ T(w_n) = \left(^{w_{n,0}}b\right)\;\, \left(^{w_{n,0}w_{n,1}}b\right) \cdots\; \, \left(^{w_{n,0}w_{n,1}\cdots w_{n,k}}b\right) \]
 where the right hand side is interpreted as a word in the alphabet 
 $^{G/B}b$ and all products in $\langle a,d \rangle$ are reduced 
to the canonical form, as elements of the generators set above. 
The $D$-action on relations 
yields the additional set of relations, for each $x\in D=\langle a, d\rangle$ 
\[ ^xT(w_n)= \left(^{xw_{n,0}}b\right)\;\, \left(^{xw_{n,0}w_{n,1}}b\right) \cdots\; \, \left(^{xw_{n,0}w_{n,1}\cdots w_{n,k}}b\right)\] 
The same procedure computes $^xT(z_n)$.
\vspace{0.1cm}
\noindent A presentation for the group $B\times B$ is now obtained by 
using  the generating set $^{G/B}b\times ^{G/B}b$ and the  following families of relations 
(coming either from relations in $B$ or from the commutativity of the two factors):
\[ (^xTw_n, 1)=1,  (1,^xTw_n)=1, (^xTz_n,1), (1,^xTz_n)=1, 
(^xb,1)(1,^yb)=(1,y^b)(^xb,1)\]  

\vspace{0.1cm}
\noindent
The next step is to obtain a presentation $P_A$  for $A$ and then 
using Hall's lemma to recover the presentation of $G$. Several remarks are in 
order. Since we seek for the finite equivalence class we can discard 
or adjoin finitely many relations at the end. When shifting from $A$ to 
$G$ we have to add the extra generators from $G/A=\langle a, c\rangle$, thus 
the generators $a$ and $c$.  We have also to add finitely many conjugation 
relations corresponding to the normality of $A$ and lifts of relations 
in $G/A$. However the previous remark enables us to ignore all these and keep 
track only of the following (four) infinite families of relations 
in $B\times B$ expressed by  
$(Tw_n, 1)=1,  (1,Tw_n)=1, (Tz_n,1)=1, (1,Tz_n)=1$.

\vspace{0.1cm}
\noindent In order to understand the isomorphism $\psi$ we have to shift 
to the presentation $P_G(a,c,d)$ of $G$.  
A natural system of generators for $A$ is given in the spirit 
of the induction lemma by the set 
$^{G/A}d=^{\langle a, c \rangle}d$ which is simply a notation for 
\[ \{^{x}d; x\in G/A =\langle a, c\rangle\}= ^{G^0/A}d \cup \;\, 
^{G^0a/A}d\]
This system of generators is convenient because $\psi$ has now a simple 
expression: 
\begin{lemma}
The isomorphism $\psi:A\to B\times B$ takes the form: 
\[ \psi_0:\{^{G^0/A}d\} \to \{1\}\times \{^{G/B}b\}, \,\,  
\psi_1:\{^{G^0a/A}d\} \to  \{^{G/B}b \}\times \{1\}\]
where 
\begin{itemize}
\item $\psi_0$ is given by: 
\[ \psi_0^{-1}(b,1)=d,\,\,  \psi_0^{-1}(^xb,1)=^{\varphi_0(x)}d \] 
where $\varphi_0:G/B=\langle a,d \rangle\to G^0/A=\langle c, aca \rangle \,$ 
is the isomorphism: 
\[ \varphi_0(d)=c,\,\,  \varphi_0(a)=aca \]
\item $\psi_1$ is given by: 
\[ \psi_1^{-1}(1,b)=^ad,\,\, \psi_0^{-1}(1,^xb)=^{\varphi_1(x)a}d \]
where  $\varphi_1:G/B=\langle a,d \rangle\to G^0/A=\langle c, aca \rangle \,$ 
is the conjugated isomorphism: 
\[ \varphi_0(d)=aca,\,\,  \varphi_0(a)=c \]
\end{itemize}
\end{lemma}
\begin{proof}
This is direct calculation. For instance $\psi(^cd)=(a,d)(1,b)(a,d)=(1,^db)$.  
\end{proof}

\vspace{0.1cm}
\noindent 
Let us transport now the relation $(1,Tw_n)=1$ from $B\times B$ to $A$. 
This relation reads: 
\[(1, ^{w_{n,0}}b)\;\, (1,^{w_{n,0}w_{n,1}}b) \cdots\; \,(1, ^{w_{n,0}w_{n,1}\cdots w_{n,k}}b) =1\]
According to the previous lemma this relation reads now in $A$ as: 
\[ \left(^{\varphi_0(w_{n,0})}d\right)\;\, \left(^{\varphi_0(w_{n,0}w_{n,1})}d\right) \cdots\; \, \left(^{\varphi_0(w_{n,0}w_{n,1}\cdots w_{n,k})}d\right) =1 \]
Further we interpret these relations in $G$ (as part of the presentation 
$P^*_G$), where we restored also the generators $a$  and $c$.  
If one writes down the terms by developing each conjugation we 
obtain: 
\[  \varphi_0(w_{n,0})d\cdot{\varphi_0(w_{n,1})}d \cdots\; \, {\varphi_0(w_{n,k})}d (\varphi_0\left(w_{n,0}w_{n,1}\cdots w_{n,k})^{-1}\right) = 1 \]
The key point  is the fact that the map $\varphi_0$ acts  
like $\sigma$ on the letters $a,d$; actually, if one extends $\varphi_0$ 
to a monoid homomorphism 
sending $b$ into $d$ we obtain $\sigma$. Thus the relation above  
is the same as: 
\[\sigma(w_n)= 1 \] 
But $\sigma(w_n)=w_{n+1}$ and thus we have no additional relation induced 
in $P^*_G$ other than those already existing in $P_G$. 

\vspace{0.1cm}
\noindent
Let us look now at the transformations of the relation 
$(1, ^xTw_n)=1$ for $x\in D$. This relation reads now in $A$ as: 
\[ \left(^{\varphi_0(xw_{n,0})}d\right)\;\, \left(^{\varphi_0(xw_{n,0}w_{n,1})}d\right) \cdots\; \, \left(^{\varphi_0(xw_{n,0}w_{n,1}\cdots w_{n,k})}d\right) =1 \]
and by developing it again in $G$: 
\[  \varphi_0(xw_{n,0})d\cdot{\varphi_0(w_{n,1})}d \cdots\; \, {\varphi_0(w_{n,k})}d (\varphi_0\left(w_{n,0}w_{n,1}\cdots w_{n,k})^{-1}\right)\varphi_0(x)^{-1} = 1 \]
This is precisely the relation: 
\[ \varphi_0(x) \sigma(w_n)\varphi_0(x)^{-1}=1 \]
which is a conjugation of the already existing relation $w_{n+1}=0$. 

\vspace{0.1cm}
\noindent
The same reasoning shows that starting from $(1,^xTz_n)$ we obtain 
in $P^*_G$ the relation $z_{n+1}=1$ (or conjugacies of it).  

\vspace{0.1cm}
\noindent Eventually we consider the relations 
$(Tw_n,1)=1$ in $B\times B$, namely: 
\[ (^{w_{n,0}}b),1)(^{w_{n,0}w_{n,1}}b, 1) \cdots\; \,(^{w_{n,0}w_{n,1}\cdots w_{n,k}}b,1) =1\]
The image of $\psi^{-1}$ of this relation in $A$ is therefore: 
\[ \left(^{\varphi_1(w_{n,0}A)}d\right)\;\, \left(^{\varphi_1(w_{n,0}w_{n,1}a)}d\right) \cdots\; \, \left(^{\varphi_1(w_{n,0}w_{n,1}\cdots w_{n,k}a)}d\right) =1 \]
But $\varphi_1(x)=a\varphi_0(x)a$ and thus this relation is the same as: 
\[ \left(^{a\varphi_0(w_{n,0})}d\right)\;\, \left(^{a\varphi_0(w_{n,0}w_{n,1})}d\right) \cdots\; \, \left(^{a\varphi_1(w_{n,0}w_{n,1}\cdots w_{n,k})}d\right) =1 \]
which, by developing all terms, yields in $G$: 
\[  a\varphi_0(w_{n,0}))d\cdot{\varphi_0(w_{n,1})}d \cdots\; \, {\varphi_0(w_{n,0}w_{n,k})}d (\varphi_0\left(w_{n,0}w_{n,1}\cdots w_{n,k})^{-1}\right)a = 1
\]

\vspace{0.1cm}
\noindent 
However this is the same as  $aw_{n+1}a=1$, which is a consequence of 
$w_{n+1}=1$. The same holds true for the relations induced by 
 $(Tz_n,1)=1$. Starting from $(^xTw_n,1)=1$ or $(^xTz_n,1)=1$
 we obtain again conjugated relations. 

\begin{lemma}
Consider two presentations of some group $G$ of the form 
$P_1=\langle S|R\rangle$ and 
$P_2=\langle S|R\cup aRa^{-1}\rangle$. 
If $P_2$ is qsf then $P_1$ is qsf.  
\end{lemma}
\begin{proof}
We can assume that $a\in S$. Every homotopy involving 
$aRa^{-1}$ can be realized using only relations from $R$. 
This proves the claim.   
\end{proof}

\vspace{0.1cm}
\noindent Now $P^*_G$ is finitely equivalent to the presentation  
consisting of $P_G$ with finitely many  additional families, each 
additional family being conjugated to the family of relations 
$\{w_n=z_n=1, n\geq 1\}$. 
If we remove the additional relations we obtain $P_G$.
The previous lemma  and Proposition \ref{finequiv} show that $P_G$ is qsf.  
This settles Proposition \ref{infgri}.

\bibliographystyle{plain}

\begin{thebibliography}{10}

\bibitem{Abe}
 H.Abels, {\em An example of a finitely presented solvable group},
Homological group theory (Proc. Sympos., Durham, 1977), pp. 205--211,
London Math. Soc. Lecture Note Ser., 36,
Cambridge Univ. Press, Cambridge-New York, 1979. 

\bibitem{AB}
J.M.Alonso and M.R.Bridson,  {\em Semihyperbolic groups}, Proc.
London Math. Soc. (3) {70}(1995), 56--114.


\bibitem{Bar}
L.Bartholdi, {\em Endomorphic presentations of branch groups},
J.Algebra  {268}(2003),  419--443.

\bibitem{BQ}
H.J.Baues and  A.Quintero,  {\em Infinite Homotopy Theory},
K-monographs in Mathematics, Kluwer Academic Publishers, 2001.



\bibitem{BF}
M.Bestvina and M.Feighn, {\em  The topology at infinity of ${\rm Out}(F\sb n)$}, 
 Invent. Math.  140(2000),  651--692.

\bibitem{Brick}
S.G.Brick, {\em Quasi-isometries and amalgamations of tame
combable groups},  Internat. J. Algebra Comput.  5(1995),
199--204.






\bibitem{BM1}
S.G.Brick and M.L.Mihalik,{\em  The QSF property for groups and
spaces},  Math. Zeitschrift   {220}(1995),  207--217.

\bibitem{BM2}
S.G.Brick and M.L.Mihalik, {\em Group extensions are
  quasi-simply-filtrated},
 Bull. Austral. Math. Soc.  50(1994), 21--27.




\bibitem{BG}
K.S.Brown and R.Geoghegan, {\em An infinite-dimensional
torsion-free ${\rm FP}\sb{\infty }$ group}, Invent. Math. 77(1984),   367--381.

\bibitem{Br}
K.S.Brown, {\em  The geometry of finitely presented infinite simple groups}, 
Algorithms and classification in combinatorial group theory (Berkeley, CA, 1989), 121--136, (G.Baumslag and C.F.Miller Ed.)
Math. Sci. Res. Inst. Publ., 23,
Springer, New York, 1992. 



\bibitem{CFP}
J.W.Cannon, W.J.Floyd and W.R.Parry, {\em  Introductory notes on
Richard Thompson's groups},
 Enseign. Math. (2)  42(1996),  215--256.




\bibitem{CM}
D.Collins and C.Miller, {\em The word problem in groups of
cohomological dimension 2},  Groups St. Andrews 1997 in Bath, I,  211--218, 
London Math. Soc. Lecture Note Ser., 260, Cambridge Univ. Press, 
Cambridge, 1999.






\bibitem{Da1}
M.W.Davis, {\em Groups generated by reflections
and aspherical manifolds not covered by Euclidian Spaces}, Ann. of
Math. 117(1983), 293--324.


\bibitem{Da2}
M.W.Davis, {\em  Exotic aspherical manifolds},
  Topology of high-dimensional manifolds, No. 1, 2 (Trieste, 2001),  371-404,
  ICTP Lect. Notes, 9, Abdus Salam Int. Cent. Theoret. Phys., Trieste, 2002.



\bibitem{DM}
M.W.Davis and J.Meier, {\em The topology at infinity of Coxeter
groups and buildings},
  Commentarii  Math. Helv. 77(2002), 746--766.



\bibitem{DLH}
P.De La Harpe, Topics in geometric group theory, Chicago Lect. Math, Chicago 
Univ. Press, 2000. 



\bibitem{DCG}
P.De La Harpe, T.Ceccherini-Silberstein and R.Grigorchuk, 
{\em Amenability and paradoxical decompositions for pseudogroups and discrete metric spaces}, Tr. Mat. Inst. Steklova  224(1999),  
Algebra. Topol. Differ. Uravn. i ikh Prilozh., 68--111;  translation in
Proc. Steklov Inst. Math. 1999, no. 1 (224), 57--97.  

\bibitem{DV}
E.Dyer and A.T.Vasquez, {\em Some small aspherical spaces}, 
J.Australian Math.Soc.  16(1973), 332--352.


\bibitem{Fa1}
D.S. Farley, {\em  Finiteness and $\rm CAT(0)$ properties of
diagram groups},   Topology  42(2003),  1065--1082.

\bibitem{Fa2}
D.S. Farley, {\em
 Actions of Picture Groups on CAT(0) Cubical Complexes},
Geometriae Dedicata 110(2005), 221--242.


\bibitem{FH}
P.Freyd and A.Heller, {\em Splitting homotopy idempotents. II}, J.
Pure Appl. Algebra {89} (1993),  93--106.








\bibitem{FG}
L.Funar and S.Gadgil, {\em  On the geometric simple connectivity
of open manifolds}, I.M.R.N.  2004, n.24, 1193--1248.


\bibitem{FLR}
L.Funar, F.Lasheras and D.Repovs, {\em  Non-compact 3-manifolds proper homotopy equivalent to geometrically simply connected polyhedra and proper 3-realizability of groups},  19p., math.GT/0709.1576.


\bibitem{FO}
L.Funar and D.E.Otera, {\em  A refinement of the simple
connectivity at infinity of groups}, Archiv Math. (Basel) 
81(2003), 360--368.






\bibitem{GS}
S.M.Gersten and J.R.Stallings, {\em Casson's idea about
$3$-manifolds whose universal cover is $\mathbb{R}^3$},
Internat. J. of Algebra Comput. 1(1991), 395--406.





\bibitem{Gri1}
R.I.Grigorchuk, {\em Degrees of growth of finitely generated groups and 
the theory of invariant means}, 
Izv. Akad. Nauk SSSR Ser. Mat. 48 (1984), no. 5, 939--985. 

\bibitem{Gri0}
R.I.Grigorchuk, {\em On a problem of M. Day on nonelementary amenable 
groups in the class of finitely presented groups},   
Mat. Zametki  60  (1996),  no. 5, 774--775;  
translation in  Math. Notes  60  (1996),  no. 5-6, 580--582. 
 

\bibitem{Gri}
R.I.Grigorchuk, {\em  An example of a finitely presented amenable
group that does not belong to the class EG}, Mat. Sb.  189(1998),
79--100; translation in  Sb. Math. 189(1998),   75--95.










\bibitem{Gr2}
M.Gromov, {\em  Asymptotic invariants of infinite groups},
Geometric group theory, Vol. 2 (Sussex, 1991), 1--295, London
Math. Soc. Lecture Note Ser. 182,  Cambridge Univ. Press,
Cambridge, 1993.



\bibitem{GS1} V.Guba and M.Sapir, {\em Diagram groups}, Mem.
Amer. Math. Soc. 130(1997),  no. 620, viii+117 pp.









\bibitem{Hig}
G.Higman, {\em A finitely generated infinite simple group},
J.London Math. Soc. 26(1951), 61--64.



\bibitem{HM1}
S.Hermiller and J.Meier, {\em
 Tame combings, almost convexity and rewriting systems for groups},
Math. Zeitschfrit  {225}(1997), 263--276.

\bibitem{HM2}
S.Hermiller and J.Meier, {\em Measuring the tameness of almost
convex groups}, Trans. Amer. Math. Soc. 353(2001),
943--962.



\bibitem{Ly}
I.G.Lys\"enok, {\em A set of defining relations for the Grigorchuk
group}, Mat. Zam. 38(1985),  503--516; transl. in Math.
Notes {38}(1985),  784--792.


\bibitem{Mi1}
M.L.Mihalik, {\em Ends of groups with integers as a quotient},
J.Pure Appl. Alg.  35(1985), 305--320.


\bibitem{Mi2}
M.L.Mihalik, {\em Solvable groups that are simply connected at
$\infty$}, Math. Zeitschfrit 195(1987),   79--87.

\bibitem{Mi3}
M.L.Mihalik, {\em Semistability of Artin and Coxeter groups}, J.Pure
Appl. Algebra 11(1996), 205--211.







\bibitem{MT2}
M.L.Mihalik and S.T.Tschantz, {\em Tame combings of groups},
  Trans. Amer. Math. Soc. 349(1997),   4251--4264.







\bibitem{O}
D.E.Otera, {\em Asymptotic topology of groups}, 
PhD thesis Univ.Paris-Sud and University of Palermo, February 2006. 



\bibitem{Po2}
  V.Poenaru, {\em Killing
handles of index one stably and $\pi _1 ^\infty $}, Duke Math. J. 
63(1991),  431--447.






\bibitem{PT}
V.Poenaru and C.Tanasi, {\em Some remarks on geometric simple
connectivity}, Acta Math. Hungarica 81(1998), 1--12.









\bibitem{St2}
J.R.Stallings, {\em Brick's quasi-simple filtrations for groups
and $3$-manifolds},   Geometric group theory, Vol. 1 (Sussex,
1991),  188--203, London Math. Soc. Lecture Note Ser., 181,
Cambridge Univ. Press, Cambridge, 1993.


\bibitem{Swe} 
E.Swenson, {\em A cut point theorem for $CAT(0)$ spaces}, 
J. Differential Geom. 53(1999), 327-358. 



\bibitem{Tu}
T.W.Tucker, {\em Non-compact $3$-manifolds and the missing
boundary problem}, Topology 73(1974), 267--273.




\bibitem{Wh2}
J.H.C.Whitehead. {\em Simple homotopy types},   Amer. J. Math.
72(1950), 1--57.


\end{thebibliography}

\end{document}